\documentclass[a4paper,11pt,reqno]{amsart}
\usepackage{amssymb,amsmath,hyperref}
\title[Appell function proofs]{Appell function proofs of recent and old mock theta function identities}

\author{Marioni Aronia}
\address{Department of Mathematics and Computer Science, Saint Petersburg State University, Saint Petersburg,  Russia, 199178}
\email{plim322pton@gmail.com}

\author{Eric T. Mortenson}
\address{Department of Mathematics and Computer Science, Saint Petersburg State University, Saint Petersburg,  Russia, 199178}
\email{etmortenson@gmail.com}

\newtheorem{theorem}{Theorem}

\newtheorem{corollary}[theorem]{Corollary}
\newtheorem{proposition}[theorem]{Proposition}
\theoremstyle{definition}

\numberwithin{theorem}{section} 
\numberwithin{equation}{section}

\makeatletter
\@namedef{subjclassname@2020}{\textup{2020} Mathematics Subject Classification}
\makeatother

\begin{document}

\date{6 May 2026}



\begin{abstract}
In this note we give new proofs of two recent mock theta function identities discovered by Garvan and Mukhopadhyay.  We also give a new proof of an old mock theta function identity of Watson.  Using the setting of Appell function properties as first introduced and developed by Hickerson and Mortenson, we demonstrate that the identities are similar to certain tenth-order and sixth-order mock theta function identities found in Ramanujan's lost notebook.  Our approach suggests more identities like those of Garvan and Mukhopadhyay.
\end{abstract}

\maketitle

 \tableofcontents

 %%%%%%%%
 %%%%%%%%
 %%%%%%%%

\section{Introduction}

In Ramanujan's last letter to Hardy, he gave a list of seventeen so-called mock theta functions.   Each function was defined by Ramanujan as a $q$-series convergent for $|q|<1$.  Although the $q$-series are not theta functions, they do have certain asymptotic properties similar to those of ordinary theta functions.  In the letter, he classifies the mock theta functions in terms of order, but the notion of order is not well-defined.

In order to state the contents of Ramanujan's last letter, we first recall some notation.  We have the $q$-Pochhammer notation
\begin{equation*}
(x)_n=(x;q)_n:=\prod_{i=0}^{n-1}(1-q^ix), \ \ (x)_{\infty}=(x;q)_{\infty}:=\prod_{i\ge 0}(1-q^ix),
\end{equation*}
as well as the condensed forms
\begin{equation*}
(x_1,x_2,\dots,x_k;q)_{n}=\prod_{i=1}^{k}(x_{i};q)_{n}, \ \ (x_1,x_2,\dots,x_n;q)_{\infty}=\prod_{i=1}^{n}(x_{i};q)_{\infty}.
\end{equation*}
We will use the following definition for the theta function
\begin{equation*}
j(z;q)=(z;q)_{\infty}(q/z;q)_{\infty}(q;q)_{\infty}=\sum_{n\in\mathbb{Z}}(-1)^{n}q^{\binom{n}{2}}z^{n},
\end{equation*}
where the last equality is the Jacobi triple product identity.  For frequently occurring types of theta functions, we will use the following notation.  Let $a,b$ be positive integers, then
\begin{equation*} 
J_{a,b}:=j(q^a;q^b),
 \ \overline{J}_{a,b}:=j(-q^a;q^b), \ {\text{and }}J_a:=J_{a,3a}=\prod_{i\ge 1}(1-q^{ai}).
\end{equation*}

In the letter we find the following three third-order mock theta functions 
\begin{equation*}
f_3(q):=\sum_{n\ge 0}\frac{q^{n^2}}{(-q)_n^2}, \ 
\psi_3(q):=\sum_{n= 1}^{\infty}\frac{q^{n^2}}{(q;q^2)_n}, \  
\chi_3(q):=\sum_{n= 0}^{\infty}\frac{q^{n^2}(-q)_n}{(-q^3;q^3)_n},
\end{equation*}
as well as several identities such as, \cite[Chapter 14]{ABV}:
{\allowdisplaybreaks \begin{gather}
f_{3}(q)+4\psi_{3}(-q)=\frac{J_{1}^3}{J_{2}^2}\label{equation:mockIdentity-f(q)psi(q)},\\
4\chi_{3}(q)-f_{3}(q)=3\frac{J_{3}^4}{J_{1}J_{6}^2},\label{equation:mockIdentity-chi(q)f(q)}\\
\chi_3(q)+\psi_3(-q)=\frac{J_{3}J_{4}^3}{J_{2}^2J_{12}}.\label{equation:mockIdentity-chi(q)psi(q)}
\end{gather}}%
For proving such mock theta function identities, the general plan is to convert the Eulerian, or $q$-hypergeometric form, of the $q$-series into building blocks, and then take advantage of building block properties.   There are many kinds of building blocks:  Appell functions, theta functions, Hecke-type double-sums, and Fourier coefficients of meromorphic Jacobi forms.  Here we will focus on Appell functions and theta functions.

We recall the definition of an Appell function
\begin{equation}
m(x,z;q):=\frac{1}{j(z;q)}
\sum_{n\in\mathbb{Z}}\frac{(-1)^{n}q^{\binom{n}{2}}z^n}{1-q^{n-1}xz}.
\label{equation:mxqz-def}
\end{equation}
Using the Watson--Whipple ${}_8\phi_{7}$ basic $q$-hypergeometric transformation:
\begin{align*}
1&+\sum_{n=1}^{\infty}\frac{(a , q\sqrt{a} ,-q\sqrt{a}, b ,c ,d ,e ,f ;q)_{n}}{(q,\sqrt{a} ,-\sqrt{a} ,{aq}/{c},{aq}/{d},{aq}/{e},{aq}/{f},{aq}/{g}; q)_{n}}\cdot \Big ( \frac{a^2q^2}{cdefg}\Big )^n \\
&=\frac{(aq,{aq}/{fg},{aq}/{ge},{aq}/{ef};q)_{\infty}}{({aq}/{e},{aq}/{f},{aq}/{g},{aq}/{efg};q)_{\infty}} \cdot \Big [ 1+\sum_{n=1}^{\infty}\frac{({aq}/{cd},e,f ,g ;q)_{n}}{(q,{efg}/{a} ,{aq}/{c} ,{aq}/{d}; q)_{n}}\cdot q^n\Big ],
\end{align*}
Watson showed (slightly rewritten) \cite{W3}
\begin{equation*}
f_{3}(q)=4m(-q,q;q^3)+\frac{J_{3,6}^2}{J_1}
\text{ and }
\psi_{3}(q)=-m(q,-q;-q^3)+\frac{qJ_{12}^3}{J_4J_{3,12}}.
\end{equation*}
By combining the two identities and using basic theta function identities, one produces
\begin{equation*}
f_{3}(q)+4\psi_{3}(-q)=4m(-q,q;q^3)-4m(-q,q;q^3)+\frac{J_{3,6}^2}{J_1}-4\frac{qJ_{12}^3}{J_4\overline{J}_{3,12}}=\frac{J_{1}^3}{J_{2}^2}.
\end{equation*}

Half a century later, additional mock theta functions and mock theta function identities were found in the lost notebook \cite{RLN}.    Here we are introduced to the mock theta conjectures.   These were ten conjectural identities for the fifth-order mock theta functions, and these conjectures were subsequently resolved by Hickerson \cite{AG, H1}.  We encounter the four tenth-order mock theta functions as well as their six very unusual identities \cite{C1, C2, C3}.   We also meet the sixth-order mock theta functions and identities, which were first discussed in \cite{AH, BC}.  The rest of our discussion will focus on identities found in the lost notebook for tenth-order and sixth-order identities.

%%%%%%%%%%%%%%%%%%%
%%%%%%%%%%%%%%%%%%%
%%%%%%%%%%%%%%%%%%%

This brings us to the tenth-order mock theta functions \cite{C1, C2, C3, RLN}
{\allowdisplaybreaks \begin{align*}
{\phi}_{10}(q)&:=\sum_{n\ge 0}\frac{q^{\binom{n+1}{2}}}{(q;q^2)_{n+1}}, \ \ {\psi}_{10}(q):=\sum_{n\ge 0}\frac{q^{\binom{n+2}{2}}}{(q;q^2)_{n+1}}, \\\ 
& \ \ \ \ \ {X}_{10}(q):=\sum_{n\ge 0}\frac{(-1)^nq^{n^2}}{(-q;q)_{2n}}, \ \  {\chi}_{10}(q):=\sum_{n\ge 0}\frac{(-1)^nq^{(n+1)^2}}{(-q;q)_{2n+1}},\notag
\end{align*}}%
which satisfy the following six slightly-rewritten identities.  We let $\omega$ be a primitive third root of unity.  The six identities from the lost notebook read  \cite{C1,C2, RLN}
{\allowdisplaybreaks \begin{align}
q^{2}\phi_{10}(q^9)-\frac{\psi_{10}(\omega q)-\psi_{10}(\omega^2 q)}{\omega - \omega^2}
&=-q\frac{J_{1,2}}{J_{3,6}}\frac{J_{3,15}J_{6}}{J_{3}},
\label{equation:tenth-id-1}\\
q^{-2}\psi_{10}(q^9)+\frac{\omega \phi_{10}(\omega q)-\omega^2\phi_{10}(\omega^2 q)}{\omega - \omega^2}
&=\frac{J_{1,2}}{J_{3,6}}\frac{J_{6,15}J_{6}}{J_{3}},
\label{equation:tenth-id-2}\\
X_{10}(q^9)-\frac{\omega \chi_{10}(\omega q)-\omega^2\chi_{10}(\omega^2 q)}{\omega - \omega^2}
&=\frac{\overline{J}_{1,4}}{\overline{J}_{3,12}}\frac{J_{18,30}J_{3}}{J_{6}},
\label{equation:tenth-id-3}\\
\chi_{10}(q^9)+q^{2}\frac{ X_{10}(\omega q)-X_{10}(\omega^2 q)}{\omega - \omega^2}&=-q^3\frac{\overline{J}_{1,4}}{\overline{J}_{3,12}}
\frac{J_{6,30}J_{3}}{J_{6}},
\label{equation:tenth-id-4}
\end{align}}%
and \cite{C3, RLN}
\begin{align}
\phi_{10}(q)-q^{-1}\psi_{10}(-q^4)+q^{-2}\chi_{10}(q^8)&=\frac{\overline{J}_{1,2}j(-q^2;-q^{10})}{J_{2,8}},
\label{equation:RLN-id-five}\\
\psi_{10}(q)+q\phi_{10}(-q^4)+X_{10}(q^8)&=\frac{\overline{J}_{1,2}j(-q^6;-q^{10})}{J_{2,8}}.
\label{equation:RLN-id-six}
\end{align}

What brought Ramanujan to these identities is unknown.  In fact,  Andrews and Berndt have stated \cite[p. 396]{ABV}

\smallskip
``{\em It is inconceivable that an identity such as (\ref{equation:RLN-id-five}) could be stumbled upon by a mindless search algorithm without any overarching theoretical insight.}''   

\smallskip
\noindent The six identities were first proved by Choi \cite{C1, C2, C3} using methods similar to those of Hickerson in his proof of the mock theta conjectures \cite{H1, H2}.  Identities (\ref{equation:tenth-id-1})--(\ref{equation:tenth-id-4}) were later given significantly shorter proofs by Zwegers \cite{Zw3}.

In \cite{Mo2018}, Mortenson gave short proofs of all six of Ramanujan's identities for the tenth-order mock theta functions by using a recent result with Hickerson on Appell function properties.  

\begin{theorem} \label{theorem:msplit-general-n} \cite[Theorem $3.5$]{HM} For generic $x,z,z'\in \mathbb{C}^*$ 
{\allowdisplaybreaks \begin{align*}
m(x,z;q)
&=\sum_{r=0}^{n-1} q^{{-\binom{r+1}{2}}} (-x)^r m\big({-}q^{{\binom{n}{2}-nr}} (-x)^n, z'; q^{n^2} \big)\\
&\qquad +z' \Theta_n^3  \sum_{r=0}^{n-1}
\frac{q^{{\binom{r}{2}}} (-xz)^r
j\big(-q^{{\binom{n}{2}+r}} (-x)^n z z';q^n\big)
j(q^{nr} z^n/z';q^{n^2})}
{j(xz;q) j(z';q^{n^2}) j\big(-q^{{\binom{n}{2}}} (-x)^n z';q^n )j(q^r z;q^n \big )}.
\end{align*}}%
\end{theorem}
\noindent The idea behind the proofs is straightforward.  First one obtains the Appell function forms of the tenth-order mock theta functions.  Next, one regroups the Appell functions by using Theorem \ref{theorem:msplit-general-n} as a guide.  Then, one replaces the Appell functions with the appropriate sums of quotients of theta functions given by Theorem \ref{theorem:msplit-general-n}.  Each of the six identities is then reduced to proving a theta function identity which can be verified through several applications of the three-term Weierstrass relation for theta functions \cite[(1)]{We}, \cite{Ko}:  For generic $a,b,c,d\in \mathbb{C}^*$
\begin{align}
j(ac;q)j(a/c;q)j(bd;q)j(b/d;q)&=j(ad;q)j(a/d;q)j(bc;q)j(b/c;q)\label{equation:3termWeier}\\
&\qquad +b/c \cdot j(ab;q)j(a/b;q)j(cd;q)j(c/d;q).\notag
\end{align}

Now we come to the sixth-order mock theta functions.  There are many identities for the sixth-order functions in the lost notebook \cite{AH, BC}.  Two of Ramanujan's sixth-order mock theta functions read \cite{AH}, \cite[Section 5]{HM}
\begin{align}
\phi(q)&:=\sum_{n\ge 0}\frac{(-1)^nq^{n^2}(q;q^2)_n}{(-q)_{2n}}=2m(q,-1;q^3),
\label{equation:6th-phi(q)}\\
\psi(q)&:=\sum_{n\ge 0}\frac{(-1)^nq^{(n+1)^2}(q;q^2)_n}{(-q)_{2n+1}}=m(1,-q;q^3).
\label{equation:6th-psi(q)}
\end{align}
They satisfy the following identities found in the lost notebook \cite[p. 135, Entry 7.4.2]{ABV}, \cite[p. 13, equations 5b, 6b]{RLN}
\begin{gather}
\phi(q^9)-\psi(q)-q^{-3}\psi(q^9)=\frac{\overline{J}_{3,12}J_{6}^2}{\overline{J}_{1,4}\overline{J}_{9,36}}
\label{equation:RLN6-A},\\
\frac{\psi(\omega q)-\psi(\omega^2q)}{(\omega-\omega^2)q}=\frac{\overline{J}_{1,4}\overline{J}_{9,36}J_{3,6}}{\overline{J}_{3,12}J_{6}}
\label{equation:RLN6-B}.
\end{gather}
Whereas the latter follows from the Appell function property \cite{HM, Zw2}
\begin{equation}
m(x,z_1;q)-m(x,z_0;q)=\frac{z_0(q)_{\infty}^3j(z_1/z_0;q)j(xz_0z_1;q)}{j(z_0;q)j(z_1;q)j(xz_0;q)j(xz_1;q)}, \label{equation:changing-z}
\end{equation}
the former is reminiscent of the six tenth-order identities (\ref{equation:tenth-id-1})-(\ref{equation:RLN-id-six}).  

Many more sixth-order mock theta functions are found in the lost notebook \cite{AH, BC}.  It was a natural to ask if they too enjoy identities similar to (\ref{equation:RLN6-A}).  Using Theorem \ref{theorem:msplit-general-n} as a guide, one is led to three new identities for the sixth-order mock theta functions:

\begin{theorem} \cite[Theorem 1.2]{Mo24C} \label{theorem:BLMSmain} The following identities for the sixth-order mock theta functions $\rho(q)$, $\sigma(q)$, $\lambda(q)$, $\mu(q)$, $\phi_{\_}(q)$, and $\psi_{\_}(q)$ are true
{\allowdisplaybreaks \begin{align}
q\rho(q)+q^3\rho(q^9)-2\sigma(q^9)
&=q\frac{J_{3,6}J_{3}^2}{J_{1,2}J_{9,18}},
\label{equation:newSixth-1}\\
q\lambda(q)+q^{3}\lambda(q^9)-2\mu(q^{9})
&=-\frac{J_{3,6}J_{6}^2}{\overline{J}_{1,4}\overline{J}_{9,36}},
\label{equation:newSixth-2}\\
\psi\_(q)+q^{-3}\psi\_(q^9)-\phi\_(q^9)
&=q\frac{\overline{J}_{3,12}J_{3}^2}{J_{1,2}J_{9,18}}.
\label{equation:newSixth-3}
\end{align}}%
\end{theorem}
\noindent The sixth-order mock theta functions $\rho(q)$, $\sigma(q)$, $\lambda(q)$, $\mu(q)$, $\phi_{\_}(q)$ are all found in the lost manuscript, so it is natural to puzzle over why identities (\ref{equation:newSixth-1})-(\ref{equation:newSixth-3}) are absent.  In \cite{Mo24C}, one also finds nineteen new tenth-order like identities for second-, sixth-, and eighth-order mock theta functions.  

As it turned out, the proofs of identities (\ref{equation:RLN6-A}), (\ref{equation:newSixth-1})-(\ref{equation:newSixth-3}) using Theorem \ref{theorem:msplit-general-n}, involved proving some rather difficult theta function identities.  One was unable to resolve them using the three-term Weierstrass identity (\ref{equation:3termWeier}); instead, one needed to appeal to methods of modularity and  a Maple software package developed by Frye and Garvan \cite{FG}.

\section{Results:  On identities of Garvan, Mukhopadhyay, and Watson}
Recently Garvan and Mukhopadhyay \cite{GM} discovered three new mock theta function identities \cite[(1.3)-(1.5)]{GM}.  This occured while they were obtaining Appell function proofs of Zwegers' modularity result for third-order mock theta functions \cite{Zw1}.  However, Garvan et al did not use Appell function properties to prove their newly discovered mock theta function identities.  Here we will use Theorem \ref{theorem:msplit-general-n} to give new proofs of two of their identities.  Our new proofs demonstrate that their identities fall within the setting of \cite{Mo2018, Mo24C, MU2024}.

This brings us to what we will prove.  Let us define $\zeta_{3}:=e^{2\pi i /3}$.   We change the notation of Garvan and Mukhopadhyay and also correct some typos.  The first identity \cite[(1.4)]{GM} has a sign mistake.  In our notation, the corrected version of the first identity reads
\begin{equation}
2q^{2}\omega_{3}(q^3)=-\frac{2i}{\sqrt{3}}+\frac{2}{3}\frac{J_{2}^4}{J_{6}J_{1}^2}
-\frac{4}{3}\frac{e^{-\pi i /3}}{(q^6;q^6)_{\infty}}\sum_{n\in\mathbb{Z}}\frac{(-1)^{n}\zeta_{3}^{2n}q^{n^2+n}}{1-\zeta_{3}q^{2n+1}}.
\label{equation:GM-1}
\end{equation}
The second identity \cite[(1.5)]{GM} has a sign mistake and an exponent mistake.  In our notation, the corrected version of the second identity reads
\begin{equation}
f_{3}(q^3)=\frac{1}{3}\frac{J_{1}^4}{J_{3}J_{2}^2}
+\frac{4}{3}\frac{1}{(q^3;q^3)_{\infty}}\sum_{n\in\mathbb{Z}}\frac{(-1)^{n}\zeta_{3}^nq^{\binom{n+1}{2}}}{1+q^{n}}.
\label{equation:GM-2}
\end{equation}

In our new proofs of identities (\ref{equation:GM-1}) and (\ref{equation:GM-2}) we resolve all of the encountered theta function identities by means of classical theta function identities such as the three-term Weierstrass relation (\ref{equation:3termWeier}).  There is no appeal to the methods of Frye and Garvan \cite{FG, Mo24C}.

We also give a new proof of an old identity of Watson \cite{W3}:
\begin{equation}
f_{3}(q^{8})+2q\omega_{3}(q)+2q^3\omega_{3}(-q^4)=\frac{J_{2}J_{4}^6}{J_{1}^2J_{8}^4}.
\end{equation}
The reader will note the similarity of Watson's identity and the identities for the sixth and tenth-order identities found in the lost notebook.
%%%%%%%
%%%%%%%
%%%%%%%

\section{Technical preliminaries}\label{section:prelim}

 We will frequently use the following product rearrangements without mention.
{\allowdisplaybreaks \begin{subequations}
\begin{gather}
\overline{J}_{0,1}=2\overline{J}_{1,4}=\frac{2J_2^2}{J_1},  \overline{J}_{1,2}=\frac{J_2^5}{J_1^2J_4^2},   J_{1,2}=\frac{J_1^2}{J_2},   \overline{J}_{1,3}=\frac{J_2J_3^2}{J_1J_6}, \notag\\
J_{1,4}=\frac{J_1J_4}{J_2},   J_{1,6}=\frac{J_1J_6^2}{J_2J_3},   \overline{J}_{1,6}=\frac{J_2^2J_3J_{12}}{J_1J_4J_6}.\notag
\end{gather}
\end{subequations}}%
Also following from our theta function definition are the general properties:
\begin{subequations}
{\allowdisplaybreaks \begin{gather}
j(q^n x;q)=(-1)^nq^{-\binom{n}{2}}x^{-n}j(x;q), \ \ n\in\mathbb{Z},\label{equation:1.8}\\
j(x;q)=j(q/x;q)\label{equation:1.7},\\
j(x;q)={J_1}j(x,qx,\dots,q^{n-1}x;q^n)/{J_n^n} \ \ {\text{if $n\ge 1$,}}\label{equation:1.10}\\
j(z;q)=\sum_{k=0}^{m-1}(-1)^k q^{\binom{k}{2}}z^k
j\big ((-1)^{m+1}q^{\binom{m}{2}+mk}z^m;q^{m^2}\big ),\label{equation:jsplit}\\
j(x^n;q^n)={J_n}j(x,\zeta_nx,\dots,\zeta_n^{n-1}x;q^n)/{J_1^n} \ \ {\text{if $n\ge 1$.}}\label{equation:1.12}
\end{gather}}%
\end{subequations}
\noindent  Here, $\zeta_n$ an $n$-th primitive root of unity.   We note two specializations of (\ref{equation:1.10}) for $n=2$ and $n=3$.  They read
\begin{gather}
j(x;q)=\frac{j(x;q^2)j(qx;q^2)J_{1}}{J_{2}^2},\label{equation:1.10n2}\\
j(x;q)=\frac{j(x;q^3)j(xq;q^3)j(xq^2;q^3)J_{1}}{J_{3}^3}.\label{equation:1.10n3}
\end{gather}
We collect additional well-known results about theta functions in terms of a proposition.
\begin{proposition}   For generic $x,y,z\in \mathbb{C}^*$ 
 \begin{subequations}
{\allowdisplaybreaks \begin{gather}
j(x;q)j(y;q)=j(-xy;q^2)j(-qx^{-1}y;q^2)-xj(-qxy;q^2)j(-x^{-1}y;q^2),\label{equation:H1Thm1.1}\\
%j(-x;q)j(y;q)-j(x;q)j(-y;q)=2xj(x^{-1}y;q^2)j(qxy;q^2),\label{equation:H1Thm1.2A}\\
%j(-x;q)j(y;q)+j(x;q)j(-y;q)=2j(xy;q^2)j(qx^{-1}y;q^2),\label{equation:H1Thm1.2B}\\
j(qx^3;q^3)+xj(q^2x^3;q^3)=j(-x;q)j(qx^2;q^2)/J_2.\label{equation:quintuple}
\end{gather}}
\end{subequations}
\end{proposition}
\noindent The last equation is the quintuple product identity.

%%%%%%%
%%%%%%%
%%%%%%%

The Appell function $m(x,z;q)$ satisfies several functional equations and identities, which we collect in the form of a proposition.  For more details, see \cite{HM, Zw2}.
\begin{proposition}  For generic $x,z\in \mathbb{C}^*$
{\allowdisplaybreaks \begin{subequations}
\begin{gather}
m(x,z;q)=m(x,qz;q),\label{equation:mxqz-fnq-z}\\
m(x,z;q)=x^{-1}m(x^{-1},z^{-1};q),\label{equation:mxqz-flip}\\
m(qx,z;q)=1-xm(x,z;q),\label{equation:mxqz-fnq-x}\\
m(x,z;q)=m(x,x^{-1}z^{-1};q).\label{equation:mxqz-flip-xz}
\end{gather}
\end{subequations}}
\end{proposition}

 We point out the $n=2$ and $n=3$ specializations of Theorem \ref{theorem:msplit-general-n} .  We have
\begin{corollary} \label{corollary:msplitn2zprime} For generic $x,z,z'\in \mathbb{C}^*$ 
{\allowdisplaybreaks \begin{align*}
m(x,z;q)&=m(-qx^2,z';q^4 )-q^{-1}xm(-q^{-1}x^2,z';q^4)\\
&\quad +\frac{z'J_{2}^3}{j(xz;q)j(z';q^4)}\Big [
\frac{j(-qx^2zz';q^2)j(z^2/z';q^{4})}{j(-qx^2z';q^2)j(z;q^2)}
-xz \frac{j(-q^2x^2zz';q^2)j(q^2z^2/z';q^{4})}{j(-qx^2z';q^2)j(qz;q^2)}\Big ].
\end{align*}}%
\end{corollary}
\begin{corollary} \label{corollary:msplitn3zprime} For generic $x,z,z'\in \mathbb{C}^*$ 
\begin{align*}
m(x,z;q)&=m\left (q^{3}x^3,z';q^{9}\right )
 -q^{-1}xm\left (x^3,z';q^{9}\right )+q^{-3}x^2m\left (q^{-3}x^3,z';q^{9}\right )\\
&\quad+\frac{z'J_3^3}{j(xz;q)j(z';q^{9})j(x^3z';q^3)}\Big [ 
\frac{1}{z}\frac{j(x^3zz';q^3)j(z^3/z';q^{9})}{j(z;q^3)} \\
&\qquad -\frac{x}{q}\frac{j(qx^3zz';q^3)j(q^{3}z^3/z';q^{9})}{j(qz;q^3)}
+\frac{x^2z}{q}\frac{j(q^2x^3zz';q^3)j(q^{6}z^3/z';q^{9})}{j(q^2z;q^3)}\Big ].
\end{align*}
\end{corollary}
%%%%%%%
%%%%%%%
%%%%%%%

\section{A new proof of Garvan and Mukhopadhyay's first identity}

We recall the identity that we want to prove
\begin{equation}
2q^{2}\omega_{3}(q^3)=-\frac{2i}{\sqrt{3}}+\frac{2}{3}\frac{J_{2}^4}{J_{6}J_{1}^2}
-\frac{4}{3}\frac{e^{-\pi i /3}}{(q^6;q^6)_{\infty}}\sum_{n\in\mathbb{Z}}\frac{(-1)^{n}\zeta_{3}^{2n}q^{n^2+n}}{1-\zeta_{3}q^{2n+1}}.
\label{equation:GMidentity1}
\end{equation}
We will need the theta-less Appell function form for the third-order function $\omega_{3}(q)$ \cite[Section 5]{HM}:
\begin{equation}
 \omega_{3}(q):=\sum_{n= 0}^{\infty}\frac{q^{2n(n+1)}}{(q;q^2)_{n+1}^2}
 =-q^{-1}m(q,q^2;q^6) -q^{-1}m(q,q^4;q^6).\label{equation:3rd-omega(q)}
\end{equation}
Using (\ref{equation:3rd-omega(q)}), the left-hand side of (\ref{equation:GMidentity1}) has the Appell function form:
\begin{equation}
2q^2\omega_{3}(q^3)=-2q^{-1}m(q^3,q^6;q^{18}) -2q^{-1}m(q^3,q^{12};q^{18}).\label{equation:alt3rdOmega}
\end{equation}

The plan is to write the sum on the right-hand side of (\ref{equation:GMidentity1}) in terms of the $m(x,z;q)$ notation.  Next we expand the Appell function by using Corollary \ref{corollary:msplitn3zprime}.  We then use Appell function properties and theta function identities to make the new right-hand side of (\ref{equation:GMidentity1}) look like the right-hand side of (\ref{equation:alt3rdOmega}).  So what we will really prove is 
\begin{equation}
\frac{4}{3}\frac{e^{-\pi i /3}}{(q^6;q^6)_{\infty}}\sum_{n\in\mathbb{Z}}\frac{(-1)^{n}\zeta_{3}^{2n}q^{n^2+n}}{1-\zeta_{3}q^{2n+1}}
=-\frac{2i}{\sqrt{3}}+\frac{2}{3}\frac{J_{2}^4}{J_{6}J_{1}^2}-2q^{2}\omega_{3}(q^3).
\label{equation:GMidentity1-alt}
\end{equation}

We begin by rewriting the sum on the right-hand side of (\ref{equation:GMidentity1}).  We make the sum look like the definition of the $m(x,z;q)$ function, see (\ref{equation:mxqz-def}).  We start with
 \begin{align*}
 \sum_{n\in\mathbb{Z}}\frac{(-1)^{n}\zeta_{3}^{2n}q^{n^2+n}}{1-\zeta_{3}q^{2n+1}}
 &= \sum_{n\in\mathbb{Z}}\frac{(-1)^{n}q^{2\binom{n}{2}}(\zeta_{3}^{2}q^{2})^n}{1-q^{2(n-1)}\zeta_{3}q^{3}}\\
  &= \sum_{n\in\mathbb{Z}}\frac{(-1)^{n}q^{2\binom{n}{2}}(\zeta_{3}^{2}q^{2})^n}
  {1-q^{2(n-1)}\zeta_{3}^{-1}q\zeta_{3}^{2}q^{2}}\\
& =j(\zeta_{3}^2q^{2};q^{2})m(\zeta_{3}^{-1}q,\zeta_{3}^2q^2;q^2).
\end{align*}
We then use the Appell function property (\ref{equation:mxqz-fnq-z}) to get
\begin{equation*}
m(\zeta_{3}^{-1}q,\zeta_{3}^2q^2;q^2)=m(\zeta_{3}^{-1}q,\zeta_{3}^2;q^2),
\end{equation*} 
and the theta function property (\ref{equation:1.7}) to obtain
\begin{equation*}
j(\zeta_{3}^2q^{2};q^{2})=-\zeta_{3}^{-2}j(\zeta_{3}^2;q^{2}).
\end{equation*} 
Assembling the two pieces, we then have
 \begin{equation}
 \sum_{n\in\mathbb{Z}}\frac{(-1)^{n}\zeta_{3}^{2n}q^{n^2+n}}{1-\zeta_{3}q^{2n+1}}
 =-\zeta_{3}^{-2}j(\zeta_{3}^2;q^{2})m(\zeta_{3}^{-1}q,\zeta_{3}^2;q^2).
 \label{equation:GMidentity1-altPre}
 \end{equation}

Now we want to use Corollary \ref{corollary:msplitn3zprime} to go from 
\begin{equation*}
m(\zeta_{3}^{-1}q,\zeta_{3}^2;q^2)
\end{equation*}
to the right-hand side of (\ref{equation:GMidentity1-alt}).  Here, we will do two expansions with Corollary \ref{corollary:msplitn3zprime}, one with $z'=q^{6}$ and the other with $z'=q^{12}$.  The choices of $z^{\prime}$ are dictated by the $z$-values found in the two Appell functions in the right-hand side of (\ref{equation:alt3rdOmega}).  We will then add the two expansions and simplify them using Appell function and theta function properties.

We first specialize Corollary \ref{corollary:msplitn3zprime} with $(x,z;q)\to (\zeta_{3}^{-1}q,\zeta_{3}^2;q^2)$
\begin{align}
m(\zeta_{3}^{-1}q,\zeta_{3}^2;q^2)&=m\left (q^{9},z';q^{18}\right )
 -q^{-1}\zeta_{3}^{-1}m\left (q^3,z';q^{18}\right )+q^{-4}\zeta_{3}^{-2}m\left (q^{-3},z';q^{18}\right )
 \label{equation:msplit3PreChoice}\\
&\quad +\frac{z'J_6^3}{j(\zeta_{3}q;q^2)j(z';q^{18})j(q^3z';q^6)}\Big [ 
\zeta_{3}\frac{j(q^3\zeta_{3}^2z';q^6)j(1/z';q^{18})}{j(\zeta_{3}^2;q^6)} \notag\\
&\qquad -\frac{\zeta_{3}^{-1}}{q}\frac{j(q^{5}\zeta_{3}^2z';q^6)j(q^{6}/z';q^{18})}{j(q^2\zeta_{3}^2;q^6)}
+\frac{j(q^{7}\zeta_{3}^2z';q^6)j(q^{12}/z';q^{18})}{j(q^4\zeta_{3}^2;q^6)}\Big ].\notag
\end{align}

For the first expansion, we further specialize $z'\to q^6$ in (\ref{equation:msplit3PreChoice}).  This yields
\begin{align}
m(\zeta_{3}^{-1}q,\zeta_{3}^2;q^2)&=m\left (q^{9},q^{6};q^{18}\right )
 -q^{-1}\zeta_{3}^{-1}m\left (q^3,q^{6};q^{18}\right )+q^{-4}\zeta_{3}^{-2}m\left (q^{-3},q^{6};q^{18}\right )
 \label{equation:firstIdentityStep1}\\
&\quad +\frac{q^{6}J_6^3}{j(\zeta_{3}q;q^2)j(q^6;q^{18})j(q^{9};q^6)}\Big [ 
\zeta_{3}\frac{j(q^{9}\zeta_{3}^2;q^6)j(q^{-6};q^{18})}{j(\zeta_{3}^2;q^6)}\notag \\
&\qquad -\frac{\zeta_{3}^{-1}}{q}\frac{j(q^{11}\zeta_{3}^2;q^6)j(1;q^{18})}{j(q^2\zeta_{3}^2;q^6)}
+\frac{j(q^{13}\zeta_{3}^2;q^6)j(q^{6};q^{18})}{j(q^4\zeta_{3}^2;q^6)}\Big ].\notag
\end{align}
We first rewrite the third Appell function.   We use Appell function property (\ref{equation:mxqz-flip}) to get
\begin{equation*}
m\left (q^{-3},q^{6};q^{18}\right )=q^3m\left (q^{3},q^{-6};q^{18}\right ),
\end{equation*}
and then (\ref{equation:mxqz-fnq-z}) to obtain
\begin{equation*}
m\left (q^{-3},q^{6};q^{18}\right )=q^3m\left (q^{3},q^{-6};q^{18}\right )=q^3m\left (q^{3},q^{12};q^{18}\right ).
\end{equation*}
Hence (\ref{equation:firstIdentityStep1}) becomes
\begin{align}
m(\zeta_{3}^{-1}q,\zeta_{3}^2;q^2)&=m\left (q^{9},q^{6};q^{18}\right )
 -q^{-1}\zeta_{3}^{-1}m\left (q^3,q^{6};q^{18}\right )+q^{-1}\zeta_{3}^{-2}m\left (q^{3},q^{12};q^{18}\right )
 \label{equation:firstIdentityStep2}\\
&\quad +\frac{q^{6}J_6^3}{j(\zeta_{3}q;q^2)j(q^6;q^{18})j(q^{9};q^6)}\Big [ 
\zeta_{3}\frac{j(q^{9}\zeta_{3}^2;q^6)j(q^{-6};q^{18})}{j(\zeta_{3}^2;q^6)} \notag\\
&\qquad -\frac{\zeta_{3}^{-1}}{q}\frac{j(q^{11}\zeta_{3}^2;q^6)j(1;q^{18})}{j(q^2\zeta_{3}^2;q^6)}
+\frac{j(q^{13}\zeta_{3}^2;q^6)j(q^{6};q^{18})}{j(q^4\zeta_{3}^2;q^6)}\Big ].\notag
\end{align}

We begin to rewrite the sum of the three quotients of theta functions in (\ref{equation:firstIdentityStep2}).    In the lead coefficient, we first focus on a theta function in the denominator and use (\ref{equation:1.8}) with $n=1$ to get
\begin{equation*}
j(q^9;q^6)=-q^{-3}j(q^3;q^6).
\end{equation*}
For the first term in brackets, we use (\ref{equation:1.8}) with $n=1$ to get
\begin{equation*}
j(q^{9}\zeta_{3}^2;q^6)=-\zeta_{3}^{-2}q^{-3}j(q^{3}\zeta_{3}^2;q^6)
\end{equation*}
and (\ref{equation:1.7}) and (\ref{equation:1.8}) with $n=1$ to obtain
\begin{equation*}
j(q^{-6};q^{18})=j(q^{24};q^{18})=-q^{-6}j(q^6;q^{24}).
\end{equation*}
The second term in brackets vanishes.  We recall that for $n\in\mathbb{Z}$, that $j(q^n;q)=0$.  Hence
\begin{equation*}
j(1;q^{18})=0.
\end{equation*}
For the third term in brackets, we use the theta function property (\ref{equation:1.8}) with $n=2$ to get
\begin{equation*}
j(q^{13}\zeta_{3}^2;q^6)=(-1)^2q^{-6\binom{2}{2}}(\zeta_{3}^{2}q)^{-2}j(q\zeta_{3}^2;q^6)=\zeta_{3}^{-1}q^{-8}j(q\zeta_{3}^2;q^6).
\end{equation*}
So the expansion (\ref{equation:firstIdentityStep2}) simplifies to 
\begin{align}
m(\zeta_{3}^{-1}q,\zeta_{3}^2;q^2)&=m\left (q^{9},q^{6};q^{18}\right )
 -q^{-1}\zeta_{3}^{-1}m\left (q^3,q^{6};q^{18}\right )+q^{-1}\zeta_{3}^{-2}m\left (q^{3},q^{12};q^{18}\right )
  \label{equation:firstIdentityStep3}\\
&\quad -\frac{\zeta_{3}^{-1}J_6^3}{j(\zeta_{3}q;q^2)j(q^6;q^{18})j(q^{3};q^6)}
\notag\\
&\qquad \times \Big [ 
\frac{j(q^{3}\zeta_{3}^2;q^6)j(q^{6};q^{18})}{j(\zeta_{3}^2;q^6)} 
 +q\frac{j(q\zeta_{3}^2;q^6)j(q^{6};q^{18})}{j(q^4\zeta_{3}^2;q^6)}\Big ].\notag
\end{align}
We see that a factor of $j(q^6;q^{18})$ cancels.  Hence (\ref{equation:firstIdentityStep3}) becomes our final form:
\begin{align}
m(\zeta_{3}^{-1}q,\zeta_{3}^2;q^2)&=m\left (q^{9},q^{6};q^{18}\right )
 -q^{-1}\zeta_{3}^{-1}m\left (q^3,q^{6};q^{18}\right )+q^{-1}\zeta_{3}^{-2}m\left (q^{3},q^{12};q^{18}\right )
   \label{equation:firstIdentityStep4}\\
&\quad -\frac{\zeta_{3}^{-1}J_6^3}{j(\zeta_{3}q;q^2)j(q^{3};q^6)}
 \times \Big [ 
\frac{j(q^{3}\zeta_{3}^2;q^6)}{j(\zeta_{3}^2;q^6)} 
 +q\frac{j(q\zeta_{3}^2;q^6)}{j(q^4\zeta_{3}^2;q^6)}\Big ].\notag
\end{align}

For the second expansion, we specialize $z'\to z^{12}$ in (\ref{equation:msplit3PreChoice}).   This produces
\begin{align}
m(\zeta_{3}^{-1}q,\zeta_{3}^2;q^2)&=m\left (q^{9},q^{12};q^{18}\right )
 -q^{-1}\zeta_{3}^{-1}m\left (q^3,q^{12};q^{18}\right )+q^{-4}\zeta_{3}^{-2}m\left (q^{-3},q^{12};q^{18}\right )
    \label{equation:secondIdentityStep1}\\
&+\frac{q^{12}J_6^3}{j(\zeta_{3}q;q^2)j(q^{12};q^{18})j(q^{15};q^6)}\Big [ 
\zeta_{3}\frac{j(q^{15}\zeta_{3}^2;q^6)j(q^{-12};q^{18})}{j(\zeta_{3}^2;q^6)} \notag\\
&\qquad -\frac{\zeta_{3}^{-1}}{q}\frac{j(q^{17}\zeta_{3}^2;q^6)j(q^{-6};q^{18})}{j(q^2\zeta_{3}^2;q^6)}
+\frac{j(q^{19}\zeta_{3}^2;q^6)j(1;q^{18})}{j(q^4\zeta_{3}^2;q^6)}\Big ].\notag
\end{align}
We then rewrite the third Appell function in (\ref{equation:secondIdentityStep1}) using (\ref{equation:mxqz-flip}) and (\ref{equation:mxqz-fnq-z}). This gives
\begin{equation*}
m\left (q^{-3},q^{12};q^{18}\right )
=q^{3}m\left (q^{3},q^{-12};q^{18}\right )=q^{3}m\left (q^{3},q^{6};q^{18}\right ).
\end{equation*}
Our expansion (\ref{equation:secondIdentityStep1}) then becomes
\begin{align}
m(\zeta_{3}^{-1}q,\zeta_{3}^2;q^2)&=m\left (q^{9},q^{12};q^{18}\right )
 -q^{-1}\zeta_{3}^{-1}m\left (q^3,q^{12};q^{18}\right )+q^{-1}\zeta_{3}^{-2}m\left (q^{3},q^{6};q^{18}\right )
  \label{equation:secondIdentityStep2}\\
&+\frac{q^{12}J_6^3}{j(\zeta_{3}q;q^2)j(q^{12};q^{18})j(q^{15};q^6)}\Big [ 
\zeta_{3}\frac{j(q^{15}\zeta_{3}^2;q^6)j(q^{-12};q^{18})}{j(\zeta_{3}^2;q^6)}\notag \\
&\qquad -\frac{\zeta_{3}^{-1}}{q}\frac{j(q^{17}\zeta_{3}^2;q^6)j(q^{-6};q^{18})}{j(q^2\zeta_{3}^2;q^6)}
+\frac{j(q^{19}\zeta_{3}^2;q^6)j(1;q^{18})}{j(q^4\zeta_{3}^2;q^6)}\Big ].\notag
\end{align}
We begin to rewrite the sum of the three quotients of theta functions in (\ref{equation:secondIdentityStep2}).    We first focus on a theta function in the lead coefficient.  We use (\ref{equation:1.8}) with $n=2$ to obtain
\begin{equation*}
j(q^{15};q^{6})=q^{-6\binom{2}{2}}\left ( q^{3}\right )^{-2}j(q^3;q^6)=q^{-12}j(q^3;q^6). 
\end{equation*}
For the first term in brackets, we use (\ref{equation:1.8}) with $n=2$ to get
\begin{equation*}
j(q^{15}\zeta_{3}^2;q^6)=q^{-6\binom{2}{2}}\left (q^3\zeta_{3}^2 \right )^{-2}j(q^{3}\zeta_{3}^2;q^6) 
=q^{-12}\zeta_{3}^{-4}j(q^{3}\zeta_{3}^2;q^6), 
\end{equation*}
and we use (\ref{equation:1.7}) and (\ref{equation:1.8}) with $n=1$ to obtain
\begin{equation*}
j(q^{-12};q^{18})=j(q^{30};q^{18})=-q^{-12}j(q^{12};q^{18}).
\end{equation*}
For the second term in brackets, we use (\ref{equation:1.8}) with $n=2$ to get
\begin{equation*}
j(q^{17}\zeta_{3}^2;q^6)=q^{-6\binom{2}{2}}\left (q^5\zeta_{3}^2 \right )^{-2}j(q^{5}\zeta_{3}^2;q^6) 
=q^{-16}\zeta_{3}^{-4}j(q^{5}\zeta_{3}^2;q^6), 
\end{equation*}
and we use (\ref{equation:1.7}) and (\ref{equation:1.8}) with $n=1$ to obtain
\begin{equation*}
j(q^{-6};q^{18})=j(q^{24};q^{18})=-q^{-6}j(q^{6};q^{18}).
\end{equation*}
The third term vanishes because $j(q^n;q)=0$ for all $n\in\mathbb{Z}$.  In other words,
\begin{equation*}
j(1;q^{18})=0.
\end{equation*}
Hence our expansion (\ref{equation:secondIdentityStep2}) becomes
\begin{align}
m(\zeta_{3}^{-1}q,\zeta_{3}^2;q^2)&=m\left (q^{9},q^{12};q^{18}\right )
 -q^{-1}\zeta_{3}^{-1}m\left (q^3,q^{12};q^{18}\right )+q^{-1}\zeta_{3}^{-2}m\left (q^{3},q^{6};q^{18}\right )
   \label{equation:secondIdentityStep3}\\
&\quad -\frac{J_6^3}{j(\zeta_{3}q;q^2)j(q^{6};q^{18})j(q^{3};q^6)}\notag \\
&\qquad \times \Big [ 
\frac{j(q^{3}\zeta_{3}^2;q^6)j(q^{12};q^{18})}{j(\zeta_{3}^2;q^6)} 
 -\zeta_{3}q\frac{j(q^{5}\zeta_{3}^2;q^6)j(q^{6};q^{18})}{j(q^2\zeta_{3}^2;q^6)}
\Big ].\notag
\end{align}
Because of (\ref{equation:1.7}) we have $j(q^{12};q^{18})=j(q^{6};q^{18})$.  Hence some theta functions cancel, and we arrive at our final form
\begin{align}
m(\zeta_{3}^{-1}q,\zeta_{3}^2;q^2)&=m\left (q^{9},q^{12};q^{18}\right )
 -q^{-1}\zeta_{3}^{-1}m\left (q^3,q^{12};q^{18}\right )+q^{-1}\zeta_{3}^{-2}m\left (q^{3},q^{6};q^{18}\right )
    \label{equation:secondIdentityStep4}\\
&\quad -\frac{J_6^3}{j(\zeta_{3}q;q^2)j(q^{3};q^6)}
 \times \Big [ 
\frac{j(q^{3}\zeta_{3}^2;q^6)}{j(\zeta_{3}^2;q^6)} 
 -\zeta_{3}q\frac{j(q^{5}\zeta_{3}^2;q^6)}{j(q^2\zeta_{3}^2;q^6)}
\Big ].\notag
\end{align}
We now add the two expansions (\ref{equation:firstIdentityStep4}) and  (\ref{equation:secondIdentityStep4}).  This brings us to
\begin{align}
2m(\zeta_{3}^{-1}q,\zeta_{3}^2;q^2)&=m\left (q^{9},q^{6};q^{18}\right )
 -q^{-1}\zeta_{3}^{-1}m\left (q^3,q^{6};q^{18}\right )+q^{-1}\zeta_{3}^{-2}m\left (q^{3},q^{12};q^{18}\right )
 \label{equation:firstIdentityPreFinal1}\\
&\quad -\frac{\zeta_{3}^{-1}J_6^3}{j(\zeta_{3}q;q^2)j(q^{3};q^6)}
 \times \Big [ 
\frac{j(q^{3}\zeta_{3}^2;q^6)}{j(\zeta_{3}^2;q^6)} 
 +q\frac{j(q\zeta_{3}^2;q^6)}{j(q^4\zeta_{3}^2;q^6)}\Big ]\notag\\
 &\quad + m\left (q^{9},q^{12};q^{18}\right )
 -q^{-1}\zeta_{3}^{-1}m\left (q^3,q^{12};q^{18}\right )+q^{-1}\zeta_{3}^{-2}m\left (q^{3},q^{6};q^{18}\right )\notag\\
&\quad -\frac{J_6^3}{j(\zeta_{3}q;q^2)j(q^{3};q^6)}
 \times \Big [ 
\frac{j(q^{3}\zeta_{3}^2;q^6)}{j(\zeta_{3}^2;q^6)} 
 -\zeta_{3}q\frac{j(q^{5}\zeta_{3}^2;q^6)}{j(q^2\zeta_{3}^2;q^6)}\Big ]. \notag
\end{align}
If we use the Appell function properties (\ref{equation:mxqz-fnq-x}), (\ref{equation:mxqz-flip}), (\ref{equation:mxqz-fnq-z}), we get
\begin{equation*}
m\left (q^{9},q^{6};q^{18}\right ) =1-q^{-9}m\left (q^{-9},q^{6};q^{18}\right ) 
=1-m\left (q^{9},q^{-6};q^{18}\right ) =1-m\left (q^{9},q^{12};q^{18}\right ). 
\end{equation*}
Hence the sum of the two expansions (\ref{equation:firstIdentityPreFinal1}) becomes
\begin{align}
2m(\zeta_{3}^{-1}q,\zeta_{3}^2;q^2)
&=1 -q^{-1}\zeta_{3}^{-1}\left ( m\left (q^3,q^{6};q^{18}\right ) +m\left (q^3,q^{12};q^{18}\right )\right ) 
 \label{equation:firstIdentityPreFinal2}\\
&\qquad  +q^{-1}\zeta_{3}^{-2}\left ( m\left (q^{3},q^{6};q^{18}\right )+m\left (q^{3},q^{12};q^{18}\right )\right )
\notag \\
&\qquad -\frac{\zeta_{3}^{-1}J_6^3}{j(\zeta_{3}q;q^2)j(q^{3};q^6)}
 \times \Big [ 
\frac{j(q^{3}\zeta_{3}^2;q^6)}{j(\zeta_{3}^2;q^6)} 
 +q\frac{j(q\zeta_{3}^2;q^6)}{j(q^4\zeta_{3}^2;q^6)}\Big ]\notag \\
&\qquad -\frac{J_6^3}{j(\zeta_{3}q;q^2)j(q^{3};q^6)}
 \times \Big [ 
\frac{j(q^{3}\zeta_{3}^2;q^6)}{j(\zeta_{3}^2;q^6)} 
 -\zeta_{3}q\frac{j(q^{5}\zeta_{3}^2;q^6)}{j(q^2\zeta_{3}^2;q^6)}\Big ]. \notag
\end{align}
We recall that $\zeta_{3}:=e^{\frac{2\pi i }{3}}$.  This allows us to write
\begin{equation*}
q^{-1}\zeta_{3}^{-2}-q^{-1}\zeta_{3}^{-1}
=q^{-1}\left (\zeta_{3} - \zeta_{3}^{-1}\right ) 
=q^{-1}\left( e^{\frac{2\pi i }{3}}-e^{-\frac{2\pi i }{3}}\right )
 =2 i q^{-1}\sin{\left (\tfrac{2\pi }{3}\right )}.
\end{equation*}
Thus our sum (\ref{equation:firstIdentityPreFinal2}), becomes
\begin{align}
2m(\zeta_{3}^{-1}q,\zeta_{3}^2;q^2)&=1+q^{-1}\sin{\left ( \tfrac{2\pi}{3}\right ) }2i\left ( m\left (q^3,q^{6};q^{18}\right )
 +m\left (q^3,q^{12};q^{18}\right )\right ) 
  \label{equation:firstIdentityPreFinal3}\\
&\qquad -\frac{\zeta_{3}^{-1}J_6^3}{j(\zeta_{3}q;q^2)j(q^{3};q^6)}
 \times \Big [ 
\frac{j(q^{3}\zeta_{3}^2;q^6)}{j(\zeta_{3}^2;q^6)} 
 +q\frac{j(q\zeta_{3}^2;q^6)}{j(q^4\zeta_{3}^2;q^6)}\Big ]\notag\\
&\qquad -\frac{J_6^3}{j(\zeta_{3}q;q^2)j(q^{3};q^6)}
 \times \Big [ 
\frac{j(q^{3}\zeta_{3}^2;q^6)}{j(\zeta_{3}^2;q^6)} 
 -\zeta_{3}q\frac{j(q^{5}\zeta_{3}^2;q^6)}{j(q^2\zeta_{3}^2;q^6)}\Big ].\notag 
\end{align}
We now use (\ref{equation:alt3rdOmega}) to rewrite the two Appell functions in terms of $\omega_{3}(q)$.  We get
\begin{align}
2m(\zeta_{3}^{-1}q,\zeta_{3}^2;q^2)&=1-i\sin{\left ( \tfrac{2\pi}{3}\right ) } 2q^{2}\omega_{3}(q^3)
  \label{equation:firstIdentityPreFinal4}\\
&\qquad -\frac{\zeta_{3}^{-1}J_6^3}{j(\zeta_{3}q;q^2)j(q^{3};q^6)}
 \times \Big [ 
\frac{j(q^{3}\zeta_{3}^2;q^6)}{j(\zeta_{3}^2;q^6)} 
 +q\frac{j(q\zeta_{3}^2;q^6)}{j(q^4\zeta_{3}^2;q^6)}\Big ]\notag\\
&\qquad -\frac{J_6^3}{j(\zeta_{3}q;q^2)j(q^{3};q^6)}
 \times \Big [ 
\frac{j(q^{3}\zeta_{3}^2;q^6)}{j(\zeta_{3}^2;q^6)} 
 -\zeta_{3}q\frac{j(q^{5}\zeta_{3}^2;q^6)}{j(q^2\zeta_{3}^2;q^6)}\Big ].\notag 
\end{align}
We rearrange the four quotients of theta functions in (\ref{equation:firstIdentityPreFinal4}) to make them easier to combine.  The goal is to reduce the four quotients of theta functions to a single quotient of theta functions.  We rearrange terms to get
\begin{align}
2m(\zeta_{3}^{-1}q,\zeta_{3}^2;q^2)&=1-i\sin{\left ( \tfrac{2\pi}{3}\right ) } 2q^{2}\omega_{3}(q^3)
  \label{equation:firstIdentityPreFinal5}\\
&\qquad -\frac{J_6^3}{j(\zeta_{3}q;q^2)j(q^{3};q^6)}
 \times \Big [ 
\zeta_{3}^{-1}\frac{j(q^{3}\zeta_{3}^2;q^6)}{j(\zeta_{3}^2;q^6)} 
+\frac{j(q^{3}\zeta_{3}^2;q^6)}{j(\zeta_{3}^2;q^6)} 
\Big ]\notag\\
&\qquad -\frac{J_6^3}{j(\zeta_{3}q;q^2)j(q^{3};q^6)}
 \times \Big [ 
 \zeta_{3}^{-1}q\frac{j(q\zeta_{3}^2;q^6)}{j(q^4\zeta_{3}^2;q^6)} -\zeta_{3}q\frac{j(q^{5}\zeta_{3}^2;q^6)}{j(q^2\zeta_{3}^2;q^6)}\Big ].\notag 
\end{align}
We combine the pair of theta functions in the second line of (\ref{equation:firstIdentityPreFinal5}).  Because $1+\zeta_{3}+\zeta_{3}^2=0$, we can write
\begin{equation*}
1+\zeta_{3}^{-1}=-\zeta_{3}.
\end{equation*}
Hence
\begin{equation}
-\frac{J_6^3}{j(\zeta_{3}q;q^2)j(q^{3};q^6)}
 \times \Big [ 
\zeta_{3}^{-1}\frac{j(q^{3}\zeta_{3}^2;q^6)}{j(\zeta_{3}^2;q^6)} 
+\frac{j(q^{3}\zeta_{3}^2;q^6)}{j(\zeta_{3}^2;q^6)} 
\Big ] 
=\frac{\zeta_{3}J_6^3}{j(\zeta_{3}q;q^2)j(q^{3};q^6)}
\frac{j(q^{3}\zeta_{3}^2;q^6)}{j(\zeta_{3}^2;q^6)}.
\label{equation:secondLineFinal}
\end{equation}
To combine the pair of theta functions in the third line of (\ref{equation:firstIdentityPreFinal5}), we will need a theta function identity.  First we combine fractions.  This brings us to
\begin{align}
 -&\frac{J_6^3}{j(\zeta_{3}q;q^2)j(q^{3};q^6)}
 \times \Big [ 
\zeta_{3}^{-1}q\frac{j(q\zeta_{3}^2;q^6)}{j(q^4\zeta_{3}^2;q^6)}- \zeta_{3}q\frac{j(q^{5}\zeta_{3}^2;q^6)}{j(q^2\zeta_{3}^2;q^6)}\Big ]
\label{equation:thirdLinePreFinal} \\
&=-\frac{\zeta_{3}^{-1}qJ_6^3}{j(\zeta_{3}q;q^2)j(q^{3};q^6)}
 \times \Big [ 
\frac{j(q^2\zeta_{3}^2;q^6)j(q\zeta_{3}^2;q^6)- \zeta_{3}^2j(q^4\zeta_{3}^2;q^6)j(q^{5}\zeta_{3}^2;q^6)}{j(q^4\zeta_{3}^2;q^6)j(q^2\zeta_{3}^2;q^6)}\Big ].\notag
\end{align}
We then use the theta function identity (\ref{equation:H1Thm1.1}) with $(x,y;q)\to (\zeta_{3}^2,-q;q^3)$.  This yields
\begin{equation*}
j(\zeta_{3}^2q;q^6)j(q^4\zeta_{3}^{-2};q^6)-\zeta_{3}^2j(q^4\zeta_{3}^2;q^6)j(\zeta_{3}^{-2}q;q^6)
=j(\zeta_{3}^2;q^3)j(-q;q^3).
\end{equation*}
We then use (\ref{equation:1.7}) to make a small adjustment in two theta functions.  This brings us to 
\begin{equation}
j(\zeta_{3}^2q;q^6)j(q^2\zeta_{3}^2;q^6)-\zeta_{3}^2j(q^4\zeta_{3}^2;q^6)j(\zeta_{3}^2q^5;q^6)
=j(\zeta_{3}^2;q^3)j(-q;q^3).
\label{equation:thirdLineIdentity}
\end{equation}
We then substitute (\ref{equation:thirdLineIdentity}) into (\ref{equation:thirdLinePreFinal}) to get
\begin{align}
 -&\frac{J_6^3}{j(\zeta_{3}q;q^2)j(q^{3};q^6)}
 \times \Big [ 
\zeta_{3}^{-1}q\frac{j(q\zeta_{3}^2;q^6)}{j(q^4\zeta_{3}^2;q^6)}- \zeta_{3}q\frac{j(q^{5}\zeta_{3}^2;q^6)}{j(q^2\zeta_{3}^2;q^6)}\Big ]
\label{equation:thirdLineFinal} \\
&=-\frac{\zeta_{3}^{-1}qJ_6^3}{j(\zeta_{3}q;q^2)j(q^{3};q^6)}
 \times \Big [ 
\frac{j(\zeta_{3}^2;q^3)j(-q;q^3)}{j(q^4\zeta_{3}^2;q^6)j(q^2\zeta_{3}^2;q^6)}\Big ].\notag 
\end{align}
We then substitute (\ref{equation:secondLineFinal}) and (\ref{equation:thirdLineFinal}) into (\ref{equation:firstIdentityPreFinal5}) to obtain
\begin{align}
2m(\zeta_{3}^{-1}q,\zeta_{3}^2;q^2)&=1-i\sin{\left ( \tfrac{2\pi}{3}\right ) } 2q^{2}\omega_{3}(q^3)
  \label{equation:firstIdentityPreFinal6}\\
&\qquad +\frac{\zeta_{3}J_6^3}{j(\zeta_{3}q;q^2)j(q^{3};q^6)}
\frac{j(q^{3}\zeta_{3}^2;q^6)}{j(\zeta_{3}^2;q^6)}\notag\\
&\qquad -\frac{\zeta_{3}^{-1}qJ_6^3}{j(\zeta_{3}q;q^2)j(q^{3};q^6)}
 \times \Big [ 
\frac{j(\zeta_{3}^2;q^3)j(-q;q^3)}{j(q^4\zeta_{3}^2;q^6)j(q^2\zeta_{3}^2;q^6)}\Big ].\notag 
\end{align}
Our last step is to combine the two quotients of theta functions in (\ref{equation:firstIdentityPreFinal6}).  For this we will need the three-term Weierstrass relation for theta functions  (\ref{equation:3termWeier}).  To save space, we introduce the definition:
\begin{equation*}
F(q):=\frac{\zeta_{3}J_6^3}{j(\zeta_{3}q;q^2)j(q^{3};q^6)}
\frac{j(q^{3}\zeta_{3}^2;q^6)}{j(\zeta_{3}^2;q^6)}\notag\\
 -\frac{\zeta_{3}^{-1}qJ_6^3}{j(\zeta_{3}q;q^2)j(q^{3};q^6)}
 \times \Big [ 
\frac{j(\zeta_{3}^2;q^3)j(-q;q^3)}{j(q^4\zeta_{3}^2;q^6)j(q^2\zeta_{3}^2;q^6)}\Big ].
\end{equation*}
We start by noting that $\zeta_{3}^{-1}=\zeta_{3}^2$ to get
\begin{align*}
F(q)=\frac{\zeta_{3}J_6^3}{j(\zeta_{3}q;q^2)j(q^{3};q^6)}
\Big [ \frac{j(q^{3}\zeta_{3}^2;q^6)}{j(\zeta_{3}^2;q^6)}
-\zeta_{3}q\frac{j(\zeta_{3}^2;q^3)j(-q;q^3)}{j(q^4\zeta_{3}^2;q^6)j(q^2\zeta_{3}^2;q^6)}
\Big ].
\end{align*}
We use (\ref{equation:1.10n2}) to factor a theta function in the second term.  This gives
\begin{align*}
F(q)&= \frac{\zeta_{3}J_6^3}{j(\zeta_{3}q;q^2)j(q^{3};q^6)}
\Big [ \frac{j(q^{3}\zeta_{3}^2;q^6)}{j(\zeta_{3}^2;q^6)}
-\zeta_{3}q\frac{j(\zeta_{3}^2;q^6)j(q^3\zeta_{3}^2;q^6)j(-q;q^3)}{j(q^4\zeta_{3}^2;q^6)j(q^2\zeta_{3}^2;q^6)}
\frac{J_{3}}{J_{6}^2}\Big ].
\end{align*}
We then pull out a common factor of $j(q^{3}\zeta_{3}^2;q^6)$ to get
\begin{align*}
F(q)&= \frac{\zeta_{3}J_6^3j(q^{3}\zeta_{3}^2;q^6)}{j(\zeta_{3}q;q^2)j(q^{3};q^6)}
\Big [ \frac{1}{j(\zeta_{3}^2;q^6)}
-\zeta_{3}q\frac{j(\zeta_{3}^2;q^6)j(-q;q^3)}{j(q^4\zeta_{3}^2;q^6)j(q^2\zeta_{3}^2;q^6)}
\frac{J_{3}}{J_{6}^2}\Big ].
\end{align*}
We begin to set up for the three-term Weierstrass relation  (\ref{equation:3termWeier}).  Product rearrangments allow us to write
\begin{equation*}
j(-q;q^3)\frac{J_{3}}{J_{6}^2}
=\frac{J_{2}J_{3}^2}{J_{1}J_{6}}\frac{J_{3}}{J_{6}^2}
=\frac{J_{2}J_{3}}{J_{1}J_{6}^2}\frac{J_{3}^2}{J_{6}}=\frac{J_{3,6}}{J_{1,6}}.
\end{equation*}
Using the product rearrangement and combining fractions brings us to
\begin{align*}
F(q)&= \frac{\zeta_{3}J_6^3j(q^{3}\zeta_{3}^2;q^6)}{j(\zeta_{3}q;q^2)j(q^{3};q^6)}
\Big [ \frac{j(q^4\zeta_{3}^2;q^6)j(q^2\zeta_{3}^2;q^6)J_{1,6}-\zeta_{3}qj(\zeta_{3}^2;q^6)^2J_{3,6}}
{j(\zeta_{3}^2;q^6)j(q^4\zeta_{3}^2;q^6)j(q^2\zeta_{3}^2;q^6)J_{1,6}}\Big ].
\end{align*}
We then use (\ref{equation:1.10n3}) to combine three theta functions in the denominator.  Factoring out the denominator then gives
\begin{align*}
F(q)&= \frac{\zeta_{3}J_6^3j(q^{3}\zeta_{3}^2;q^6)}{j(\zeta_{3}q;q^2)j(q^{3};q^6)}\frac{J_{2}}{J_{6}^3}\frac{1}{J_{1,6}}\frac{1}{j(\zeta_{3}^2;q^2)}
\Big [ j(q^4\zeta_{3}^2;q^6)j(q^2\zeta_{3}^2;q^6)J_{1,6}-\zeta_{3}qj(\zeta_{3}^2;q^6)^2J_{3,6}
\Big ].
\end{align*}
We now use the three-term Weierstrass relation  (\ref{equation:3termWeier}) to rewrite the sum of theta functions inside the brackets.  In (\ref{equation:3termWeier}) we set $(a,b,c,d;q)\to (\zeta_{3}^2q^3,q^3,q^2,q;q^6)$.  This gives
\begin{align*}
j(q^5\zeta_{3}^2;q^6)&j(q\zeta_{3}^2;q^6)j(q^4;q^6)j(q^2;q^6)\\
&=j(\zeta_{3}^2q^4;q^6)j(\zeta_{3}^2q^2;q^6)j(q^5;q^6)j(q;q^6)
+qj(q^6\zeta_{3}^2;q^6)j(\zeta_{3}^2;q^6)j(q^3;q^6)j(q;q^6).
\end{align*}
From (\ref{equation:1.7}) we know
\begin{equation*}
j(q^4;q^6)=j(q^2;q^6)=J_{2}, \ j(q^5;q^6)=j(q;q^6)=J_{1,6},
\end{equation*}
and from (\ref{equation:1.8}) with $n=1$, we have
\begin{equation*}
j(q^6\zeta_{3}^2;q^6)=-\zeta_{3}^{-2}j(\zeta_{3}^2;q^6)=-\zeta_{3}j(\zeta_{3}^2;q^6).
\end{equation*}
Hence
\begin{equation*}
j(q^5\zeta_{3}^2;q^6)j(q\zeta_{3}^2;q^6)J_{2}^2
=j(\zeta_{3}^2q^4;q^6)j(\zeta_{3}^2q^2;q^6)J_{1,6}^2
-\zeta_{3}qj(\zeta_{3}^2;q^6)j(\zeta_{3}^2;q^6)J_{3,6}J_{1,6}.
\end{equation*}
Thus
\begin{align*}
F(q)&= \frac{\zeta_{3}J_6^3j(q^{3}\zeta_{3}^2;q^6)}{j(\zeta_{3}q;q^2)j(q^{3};q^6)}\frac{J_{2}}{J_{6}^3}\frac{1}{J_{1,6}}\frac{1}{j(\zeta_{3}^2;q^2)}
\Big [ \frac{j(q^5\zeta_{3}^2;q^6)j(q\zeta_{3}^2;q^6)J_{2}^2}{J_{1,6}}
\Big ].
\end{align*}
We combine three theta functions in the numerator using (\ref{equation:1.10n3}).  This gives
\begin{equation*}
F(q)= \frac{\zeta_{3}J_6^3}{j(q^{3};q^6)}\frac{1}{J_{1,6}}\frac{1}{j(\zeta_{3}^2;q^2)} \frac{J_{2}^2}{J_{1,6}}.
\end{equation*}
Product rearrangements yield
\begin{equation*}
F(q)= \frac{\zeta_{3}}{j(\zeta_{3}^2;q^2)} \frac{J_{2}^4}{J_{1}^2}.
\end{equation*}
Substituting the value of $F(q)$ into (\ref{equation:firstIdentityPreFinal6}) and evaluating the sine gives
\begin{equation}
2m(\zeta_{3}^{-1}q,\zeta_{3}^2;q^2)=1-i\sqrt{3}q^{2}\omega_{3}(q^3)
 + \frac{\zeta_{3}}{j(\zeta_{3}^2;q^2)} \frac{J_{2}^4}{J_{1}^2}. \label{equation:firstIdentityPreFinal7}
\end{equation}
Using the facts that 
\begin{equation*}
j(\zeta_{3}^2;q^{2})=(1-\zeta_{3}^2)J_{6} \ \textup{and} \ \zeta_{3}-\zeta_{3}^{-1}=2i\sin{\left ( \tfrac{2\pi }{3}\right ) }
\end{equation*}
gives
\begin{equation*}
 \frac{\zeta_{3}}{j(\zeta_{3}^2;q^2)} =\frac{\zeta_{3}}{(1-\zeta_{3}^2)J_{6}}
 =-\frac{1}{(\zeta_{3}-\zeta_{3}^{-1})J_{6}}=-\frac{1}{i\sqrt{3}J_{6}}.
\end{equation*}
Hence (\ref{equation:firstIdentityPreFinal7}) becomes
\begin{equation}
2m(\zeta_{3}^{-1}q,\zeta_{3}^2;q^2)=1-i\sqrt{3}q^{2}\omega_{3}(q^3)
 - \frac{1}{i\sqrt{3}} \frac{J_{2}^4}{J_{1}^2J_{6}}. \label{equation:firstIdentityPreFinal8}
\end{equation}

Recalling (\ref{equation:GMidentity1-altPre}), we get
\begin{align*}
\frac{4}{3}\frac{e^{-\pi i /3}}{(q^6;q^6)_{\infty}}\sum_{n\in\mathbb{Z}}\frac{(-1)^{n}\zeta_{3}^{2n}q^{n^2+n}}{1-\zeta_{3}q^{2n+1}}
=
-\frac{4}{3}\frac{e^{-\pi i /3}}{(q^6;q^6)_{\infty}}
\zeta_{3}^{-2}j(\zeta_{3}^2;q^{2})m(\zeta_{3}^{-1}q,\zeta_{3}^2;q^2)
\end{align*}
Using the fact that 
\begin{equation*}
j(\zeta_{3}^2;q^{2})=(1-\zeta_{3}^2)J_{6}
\end{equation*}
gives
\begin{align*}
\frac{4}{3}\frac{e^{-\pi i /3}}{(q^6;q^6)_{\infty}}\sum_{n\in\mathbb{Z}}\frac{(-1)^{n}\zeta_{3}^{2n}q^{n^2+n}}{1-\zeta_{3}q^{2n+1}}
&=
-\frac{4}{3}e^{-\pi i /3}
\zeta_{3}^{-2}(1-\zeta_{3}^2)m(\zeta_{3}^{-1}q,\zeta_{3}^2;q^2)\\
&=\frac{4}{3}e^{-\pi i /3}
\zeta_{3}^{-1}(\zeta_{3}-\zeta_{3}^{-1})m(\zeta_{3}^{-1}q,\zeta_{3}^2;q^2)\\
&=-\frac{4}{3}2i\sin{\left ( \tfrac{2\pi }{3}\right) }m(\zeta_{3}^{-1}q,\zeta_{3}^2;q^2)\\
&=-\frac{4i}{\sqrt{3}}m(\zeta_{3}^{-1}q,\zeta_{3}^2;q^2).
\end{align*}
Using (\ref{equation:firstIdentityPreFinal8}) gives
\begin{align*}
\frac{4}{3}\frac{e^{-\pi i /3}}{(q^6;q^6)_{\infty}}\sum_{n\in\mathbb{Z}}\frac{(-1)^{n}\zeta_{3}^{2n}q^{n^2+n}}{1-\zeta_{3}q^{2n+1}}
&=-\frac{2i}{\sqrt{3}}\left  ( 1-i\sqrt{3}q^{2}\omega_{3}(q^2)
 - \frac{1}{i\sqrt{3}} \frac{J_{2}^4}{J_{1}^2J_{6}}\right )\\
&=-\frac{2i}{\sqrt{3}}-2q^{2}\omega_{3}(q^3)
 + \frac{2}{3} \frac{J_{2}^4}{J_{1}^2J_{6}}.
\end{align*}
%%%%%%%%
%%%%%%%%
%%%%%%%%

\section{A new proof of Garvan and Mukhopadhyay's second identity}
We recall the identity that we want to prove
\begin{equation}
f_{3}(q^3)=\frac{1}{3}\frac{J_{1}^4}{J_{3}J_{2}^2}
+\frac{4}{3}\frac{1}{(q^3;q^3)_{\infty}}\sum_{n\in\mathbb{Z}}\frac{(-1)^{n}\zeta_{3}^nq^{\binom{n+1}{2}}}{1+q^{n}}.
\label{equation:GMidentity2}
\end{equation}
We will need the theta-less Appell function form for the third-order function $f_3(q)$ \cite[Section 5]{HM}:
\begin{equation}
f_{3}(q):=\sum_{n\ge 0}\frac{q^{n^2}}{(-q)_n^2}=2m(-q,q;q^3)+2m(-q,q^2;q^3)\label{equation:3rd-f(q)}.
\end{equation}
Using (\ref{equation:3rd-f(q)}), the left-hand side of (\ref{equation:GMidentity2}) has the Appell function form:
\begin{equation}
f_{3}(q^3)=2m(-q^3,q^3;q^9)+2m(-q^3,q^6;q^9)\label{equation:alt3rdf}.
\end{equation}
As in the previous section, the plan is to express the sum on the right-hand side of  (\ref{equation:GMidentity2}) in terms of the $m(x,z;q)$ function.  After that, we expand the Appell function by using Corollary \ref{corollary:msplitn3zprime}.  Next, we use Appell function properties and theta function identities to make the new right-hand side of (\ref{equation:GMidentity2}) look like the right-hand side of (\ref{equation:alt3rdf}).  So our aim is actually to prove is 
\begin{equation}
\frac{4}{3}\frac{1}{(q^3;q^3)_{\infty}}\sum_{n\in\mathbb{Z}}\frac{(-1)^{n}\zeta_{3}^nq^{\binom{n+1}{2}}}{1+q^{n}}
=f_{3}(q^3)-\frac{1}{3}\frac{J_{1}^4}{J_{3}J_{2}^2}.
\label{equation:GMidentity2-alt}
\end{equation}

We proceed to express the sum on the right-hand side of  (\ref{equation:GMidentity2}) in terms of the definition of the $m(x,z;q)$ function (\ref{equation:mxqz-def}).  We begin with
\begin{align*}
\sum_{n\in\mathbb{Z}}\frac{(-1)^{n}\zeta_{3}^nq^{\binom{n+1}{2}}}{1+q^{n}}
&=\sum_{n\in\mathbb{Z}}\frac{(-1)^{n}(\zeta_{3}q)^nq^{\binom{n}{2}}}{1-q^{n-1}(-q)}\\
&=\sum_{n\in\mathbb{Z}}\frac{(-1)^{n}(\zeta_{3}q)^nq^{\binom{n}{2}}}{1-q^{n-1}(-\zeta_{3}^{-1})(\zeta_{3}q)}\\
&=j(\zeta_{3}q;q)m(-\zeta_{3}^{-1},\zeta_{3}q;q).
\end{align*}
We use the theta function property (\ref{equation:1.8}) with $n=1$ and the Appell function property (\ref{equation:mxqz-fnq-z}) to get
\begin{equation*}
\sum_{n\in\mathbb{Z}}\frac{(-1)^{n}\zeta_{3}^nq^{\binom{n+1}{2}}}{1+q^{n}}=-\zeta_{3}^{-1}j(\zeta_{3};q)m(-\zeta_{3}^{-1},\zeta_{3};q).
\end{equation*}
Using the Appell function properties (\ref{equation:mxqz-flip}) and then (\ref{equation:mxqz-flip-xz}) we get
\begin{equation*}
-\zeta_{3}^{-1}m(-\zeta_{3}^{-1},\zeta_{3};q)
=m(-\zeta_{3},\zeta_{3}^{-1};q)=m(-\zeta_{3},-1;q).
\end{equation*}
Finally we note that 
\begin{equation*}
j(\zeta_{3};q)=(1-\zeta_{3})(q^3;q^3)_{\infty}.
\end{equation*}
We thus arrive at
\begin{equation}
\frac{1}{(q^3;q^3)_{\infty}}\sum_{n\in\mathbb{Z}}\frac{(-1)^{n}\zeta_{3}^nq^{\binom{n+1}{2}}}{1+q^{n}}
=(1-\zeta_{3})m(-\zeta_{3},-1;q). \label{equation:GMidentity2-altPre}
\end{equation}

Now we will use Corollary \ref{corollary:msplitn3zprime} to go from 
\begin{equation*}
m(-\zeta_{3},-1;q)
\end{equation*}
to the right-hand side of (\ref{equation:GMidentity2-alt}).  Here, we will do two expansions with Corollary \ref{corollary:msplitn3zprime}, one with $z'=q^{3}$ and the other with $z'=q^{6}$.  The two choices of $z^{\prime}$ are dictated by the $z$-values found in the two Appell functions in the right-hand side of (\ref{equation:alt3rdf}).  We will then add the two expansions and simplify them using Appell function and theta function properties.

Corollary \ref{corollary:msplitn3zprime} with $(x,z;q)\to (-\zeta_{3},-1;q)$ gives us
\begin{align}
m(-\zeta_{3},-1;q)&=m\left (-q^{3},z';q^{9}\right )
 +q^{-1}\zeta_{3}m\left (-1,z';q^{9}\right )+q^{-3}\zeta_{3}^2m\left (-q^{-3},z';q^{9}\right )
  \label{equation:msplit3PreChoice-GM2}\\
&\quad+\frac{z'J_3^3}{j(\zeta_{3};q)j(z';q^{9})j(-z';q^3)}\Big [ 
-\frac{j(z';q^3)j(-1/z';q^{9})}{j(-1;q^3)}\notag \\
&\qquad -\frac{\zeta_{3}}{q}\frac{j(qz';q^3)j(-q^{3}/z';q^{9})}{j(-q;q^3)}
-\frac{\zeta_{3}^2}{q}\frac{j(q^2z';q^3)j(-q^{6}/z';q^{9})}{j(-q^2;q^3)}\Big ].\notag
\end{align}
For the first expansion, we specialize (\ref{equation:msplit3PreChoice-GM2}) with $z^{\prime}=q^{3}$.  This yields
\begin{align*}
m(-\zeta_{3},-1;q)&=m\left (-q^{3},q^3;q^{9}\right )
 +q^{-1}\zeta_{3}m\left (-1,q^3;q^{9}\right )+q^{-3}\zeta_{3}^2m\left (-q^{-3},q^3;q^{9}\right )\\
&\quad+\frac{q^3J_3^3}{j(\zeta_{3};q)j(q^3;q^{9})j(-q^3;q^3)}\Big [ 
-\frac{j(q^3;q^3)j(-1/q^3;q^{9})}{j(-1;q^3)} \\
&\qquad +\frac{\zeta_{3}}{q}\frac{j(q^{4};q^3)j(-1;q^{9})}{j(-q;q^3)}
-\frac{\zeta_{3}^2}{q}\frac{j(q^{5};q^3)j(-q^{3};q^{9})}{j(-q^2;q^3)}\Big ].
\end{align*}
We rewrite the third Appell function using (\ref{equation:mxqz-flip}) and (\ref{equation:mxqz-fnq-z}) to find
\begin{equation*}
q^{-3}m\left (-q^{-3},q^3;q^{9}\right )=-m\left (-q^{3},q^{-3};q^{9}\right )=-m\left (-q^{3},q^{6};q^{9}\right ).
\end{equation*}
We also note that $j(q^3;q^3)=0$ because $j(q^n;q)=0$ for $n\in\mathbb{Z}$.   Hence
\begin{align*}
m(-\zeta_{3},-1;q)&=m\left (-q^{3},q^3;q^{9}\right )
 +q^{-1}\zeta_{3}m\left (-1,q^3;q^{9}\right )-\zeta_{3}^2m\left (-q^{3},q^{6};q^{9}\right )\\
&\quad+\frac{q^3J_3^3}{j(\zeta_{3};q)j(q^3;q^{9})j(-q^3;q^3)}\\
&\qquad \times \Big [  \frac{\zeta_{3}}{q}\frac{j(q^{4};q^3)j(-1;q^{9})}{j(-q;q^3)}
-\frac{\zeta_{3}^2}{q}\frac{j(q^{5};q^3)j(-q^{3};q^{9})}{j(-q^2;q^3)}\Big ].
\end{align*}
Using the theta function property (\ref{equation:1.7}) with $n=1$ yields
\begin{align*}
m(-\zeta_{3},-1;q)&=m\left (-q^{3},q^3;q^{9}\right )
 +q^{-1}\zeta_{3}m\left (-1,q^3;q^{9}\right )-\zeta_{3}^2m\left (-q^{3},q^{6};q^{9}\right )\\
&\quad+\frac{q^3J_3^3}{j(\zeta_{3};q)j(q^3;q^{9})j(-q^3;q^3)}\\
&\qquad \times \Big [  -q^{-1}\frac{\zeta_{3}}{q}\frac{j(q;q^3)j(-1;q^{9})}{j(-q;q^3)}
+q^{-2}\frac{\zeta_{3}^2}{q}\frac{j(q^{2};q^3)j(-q^{3};q^{9})}{j(-q^2;q^3)}\Big ].
\end{align*}
From (\ref{equation:1.8}), we know that
\begin{equation*}
j(q^2;q^3)=j(q;q^3), \ j(-q^2;q^3)=j(-q;q^3),
\end{equation*} 
so we can factor out common terms and obtain
\begin{align}
m(-\zeta_{3},-1;q)&=m\left (-q^{3},q^3;q^{9}\right )
 +q^{-1}\zeta_{3}m\left (-1,q^3;q^{9}\right )-\zeta_{3}^2m\left (-q^{3},q^{6};q^{9}\right )
 \label{equation:GM2-spec1}\\
&\quad+\frac{J_3^3j(q;q^3)}{j(\zeta_{3};q)j(q^3;q^{9})j(-1;q^3)j(-q;q^3)}\notag\\
&\qquad \times \Big [  -\zeta_{3}qj(-1;q^{9})
+\zeta_{3}^2j(-q^{3};q^{9})\Big ].\notag
\end{align}

For the second expansion, we specialize (\ref{equation:msplit3PreChoice-GM2}) with $z^{\prime}=q^{6}$.  This yields
\begin{align*}
m(-\zeta_{3},-1;q)&=m\left (-q^{3},q^{6};q^{9}\right )
 +q^{-1}\zeta_{3}m\left (-1,q^{6};q^{9}\right )+q^{-3}\zeta_{3}^2m\left (-q^{-3},q^{6};q^{9}\right )\\
&\quad+\frac{q^{6}J_3^3}{j(\zeta_{3};q)j(q^{6};q^{9})j(-q^{6};q^3)}\Big [ 
-\frac{j(q^{6};q^3)j(-1/q^{6};q^{9})}{j(-1;q^3)} \\
&\qquad +\frac{\zeta_{3}}{q}\frac{j(q^{7};q^3)j(-q^{-3};q^{9})}{j(-q;q^3)}
-\frac{\zeta_{3}^2}{q}\frac{j(q^{8};q^3)j(-1;q^{9})}{j(-q^2;q^3)}\Big ].
\end{align*}
We rewrite the third Appell function using (\ref{equation:mxqz-flip}) and (\ref{equation:mxqz-fnq-z}).  This gives
\begin{equation*}
q^{-3}m\left (-q^{-3},q^6;q^{9}\right )=-m\left (-q^{3},q^{-6};q^{9}\right )=-m\left (-q^{3},q^{3};q^{9}\right ).
\end{equation*}
We also note that $j(q^6;q^3)=0$ because $j(q^n;q)=0$ for $n\in\mathbb{Z}$.   Hence
\begin{align*}
m(-\zeta_{3},-1;q)&=m\left (-q^{3},q^{6};q^{9}\right )
 +q^{-1}\zeta_{3}m\left (-1,q^{6};q^{9}\right )-\zeta_{3}^2m\left (-q^{3},q^{3};q^{9}\right )\\
&\quad+\frac{q^{6}J_3^3}{j(\zeta_{3};q)j(q^{6};q^{9})j(-q^{6};q^3)}  \\
&\qquad \times \Big [ \frac{\zeta_{3}}{q}\frac{j(q^{7};q^3)j(-q^{-3};q^{9})}{j(-q;q^3)}
-\frac{\zeta_{3}^2}{q}\frac{j(q^{8};q^3)j(-1;q^{9})}{j(-q^2;q^3)}\Big ].
\end{align*}
From (\ref{equation:1.7}), we have
\begin{equation*}
j(q^{-3};q^9)=j(q^{12};q^{9}), \ j(q^{6};q^{9})=j(q^{3};q^{9}).
\end{equation*}
Thus
\begin{align*}
m(-\zeta_{3},-1;q)&=m\left (-q^{3},q^{6};q^{9}\right )
 +q^{-1}\zeta_{3}m\left (-1,q^{6};q^{9}\right )-\zeta_{3}^2m\left (-q^{3},q^{3};q^{9}\right )\\
&\quad+\frac{q^{6}J_3^3}{j(\zeta_{3};q)j(q^{3};q^{9})j(-q^{6};q^3)}  \\
&\qquad \times \Big [ \frac{\zeta_{3}}{q}\frac{j(q^{7};q^3)j(-q^{12};q^{9})}{j(-q;q^3)}
-\frac{\zeta_{3}^2}{q}\frac{j(q^{8};q^3)j(-1;q^{9})}{j(-q;q^3)}\Big ].
\end{align*}
We then use (\ref{equation:1.8}) to rewrite the theta functions
\begin{equation*}
j(-q^6;q^3), \ j(q^7;q^3), \ j(-q^{12};q^{9}), \ j(q^{8};q^{3}).
\end{equation*}
This gives
\begin{align*}
m(-\zeta_{3},-1;q)&=m\left (-q^{3},q^{6};q^{9}\right )
 +q^{-1}\zeta_{3}m\left (-1,q^{6};q^{9}\right )-\zeta_{3}^2m\left (-q^{3},q^{3};q^{9}\right )\\
&\quad+\frac{q^{9}J_3^3}{j(\zeta_{3};q)j(q^{3};q^{9})j(-1;q^3)}  \\
&\qquad \times \Big [ q^{-9}\zeta_{3}\frac{j(q;q^3)j(-q^{3};q^{9})}{j(-q;q^3)}
-q^{-8}\zeta_{3}^2\frac{j(q^{2};q^3)j(-1;q^{9})}{j(-q;q^3)}\Big ].
\end{align*}
From (\ref{equation:1.7}), we know
\begin{equation*}
j(q^2;q^3)=j(q;q^3),
\end{equation*}
hence we can pull out common factors to get
\begin{align}
m(-\zeta_{3},-1;q)&=m\left (-q^{3},q^{6};q^{9}\right )
 +q^{-1}\zeta_{3}m\left (-1,q^{6};q^{9}\right )-\zeta_{3}^2m\left (-q^{3},q^{3};q^{9}\right )
  \label{equation:GM2-spec2}\\
&\quad+\frac{J_3^3j(q;q^3)}{j(\zeta_{3};q)j(q^{3};q^{9})j(-1;q^3)j(-q;q^3)}\notag  \\
&\qquad \times \Big [ \zeta_{3}j(-q^{3};q^{9})
-q\zeta_{3}^2j(-1;q^{9})\Big ].\notag
\end{align}
Adding the (\ref{equation:GM2-spec1}) and (\ref{equation:GM2-spec2}) brings us to
\begin{align*}
2m(-\zeta_{3},-1;q)&=\left ( 1 - \zeta_{3}^2 \right ) m\left (-q^{3},q^3;q^{9}\right )
 +q^{-1}\zeta_{3}\left ( m\left (-1,q^3;q^{9}\right )+m\left (-1,q^{6};q^{9}\right )\right )\\
 &\qquad  +\left ( 1-\zeta_{3}^2\right ) m\left (-q^{3},q^{6};q^{9}\right )\\
&\quad+\frac{J_3^3j(q;q^3)}{j(\zeta_{3};q)j(q^3;q^{9})j(-1;q^3)j(-q;q^3)}\\
&\qquad \times \Big [ -\left (  \zeta_{3}+\zeta_{3}^2\right ) qj(-1;q^{9})
+\left ( \zeta_{3}^2+\zeta_{3}\right ) j(-q^{3};q^{9})\Big ].
\end{align*}
Appell function properties (\ref{equation:mxqz-flip}) and (\ref{equation:mxqz-fnq-z}) give
\begin{equation*}
m(-1,q^{6};q^{9})=-m(-1,q^{-6};q^{9})=-m(-1,q^{3};q^{9}).
\end{equation*}
We also have that
\begin{equation*}
\zeta_{3}+\zeta_{3}^2=-1.
\end{equation*}
Hence
\begin{align*}
2m(-\zeta_{3},-1;q)&=
\left ( 1 - \zeta_{3}^2 \right ) \left ( m\left (-q^{3},q^3;q^{9}\right ) -m\left (-q^{3},q^6;q^{9} \right )\right )  \\
&\quad- \frac{J_3^3j(q;q^3)}{j(\zeta_{3};q)j(q^3;q^{9})j(-1;q^3)j(-q;q^3)}
 \times \Big [
 j(-q^{3};q^{9}) - qj(-1;q^{9})\Big ].
\end{align*}

We recall the quintuple product identity (\ref{equation:quintuple}) and set $(x;q)\to(-q;q^3)$.  This gives
\begin{equation*}
j(-q^6;q^9)-qj(-q^{9};q^9)=j(q;q^3)j(q^{5};q^6)/J_{6}.
\end{equation*}
Using (\ref{equation:1.7}), we can adjust the two theta functions: 
\begin{equation*}
j(-q^3;q^9)-qj(-1;q^9)=j(q;q^3)j(q;q^6)/J_{6}.
\end{equation*}
and as a result
\begin{align*}
2m(-\zeta_{3},-1;q)&=\left ( 1 - \zeta_{3}^2 \right ) 
\left ( m\left (-q^{3},q^3;q^{9}\right ) - m\left (-q^{3},q^{6};q^{9}\right ) \right ) \\
&\qquad -\frac{J_3^3j(q;q^3)}{j(\zeta_{3};q)j(q^3;q^{9})j(-1;q^3)j(-q;q^3)}
\times \frac{j(q;q^3)j(q;q^{6})}{J_{6}}.
\end{align*}
We use the Appell function form of $f_{3}(q^3)$ in (\ref{equation:alt3rdf}) to obtain
\begin{align*}
2m(-\zeta_{3},-1;q)&= \left ( 1 - \zeta_{3}^2 \right ) \frac{1}{2} f_{3}(q^3)\\
&\quad -\frac{J_3^3j(q;q^3)}{j(\zeta_{3};q)j(q^3;q^{9})j(-1;q^3)j(-q;q^3)}
\times \frac{j(q;q^3)j(q;q^{6})}{J_{6}}.
\end{align*}
It is straightforward to see that
\begin{equation*}
j(\zeta_{3};q)=(1-\zeta_{3})J_{3}.
\end{equation*}
Elementary product rearrangements then give
\begin{align*}
2m(-\zeta_{3},-1;q)&= \left ( 1 - \zeta_{3}^2 \right ) \frac{1}{2} f_{3}(q^3)
-\frac{1}{(1-\zeta_{3})}\frac{1}{2}\frac{J_{1}^4}{J_{3}J_{2}^2}.
\end{align*}
Substituting into (\ref{equation:GMidentity2-altPre}) gives
\begin{align*}
\frac{4}{3}\frac{1}{(q^3;q^3)_{\infty}}\sum_{n\in\mathbb{Z}}\frac{(-1)^{n}\zeta_{3}^nq^{\binom{n+1}{2}}}{1+q^{n}}
&=\frac{2}{3}(1-\zeta_{3}) \left ( \left ( 1 - \zeta_{3}^2 \right ) \frac{1}{2} f_{3}(q^3)
-\frac{1}{(1-\zeta_{3})}\frac{1}{2}\frac{J_{1}^4}{J_{3}J_{2}^2}\right ). 
\end{align*}
We note that
\begin{equation*}
\left ( 1 - \zeta_{3} \right ) \left ( 1 - \zeta_{3}^2 \right )=1-\zeta_{3}-\zeta_{3}^2+1=3.
\end{equation*}
Hence
\begin{align*}
\frac{4}{3}\frac{1}{(q^3;q^3)_{\infty}}\sum_{n\in\mathbb{Z}}\frac{(-1)^{n}\zeta_{3}^nq^{\binom{n+1}{2}}}{1+q^{n}}
&=f_{3}(q^3) -\frac{1}{3}\frac{J_{1}^4}{J_{3}J_{2}^2}. 
\end{align*}

%%%%%%%%
%%%%%%%%
%%%%%%%%

\section{A new proof of Watson's identity}
We want to prove
\begin{equation}
f_{3}(q^{8})+2q\omega_{3}(q)+2q^3\omega_{3}(-q^4)=\frac{J_{2}J_{4}^6}{J_{1}^2J_{8}^4}.
\end{equation}
However, we will rewrite the identity as
\begin{equation}
-2q\omega_{3}(q)=f_{3}(q^{8})+2q^3\omega_{3}(-q^4)-\frac{J_{2}J_{4}^6}{J_{1}^2J_{8}^4}.
\label{equation:firstIdWatson}
\end{equation}
We recall the theta-less Appell function expansions for $f_{3}(q)$ and $\omega_{3}(q)$ as found in \cite[Section 5]{HM}:
\begin{gather}
f_{3}(q)=2m(-q,q;q^3)+2m(-q,q^2;q^3),\label{equation:thetaLessf3}\\
q\omega_{3}(q)=-m(q,q^2;q^6)-m(q,q^4;q^6)\label{equation:thetaLessOmega3}.
\end{gather}
The plan is to apply Corollary \ref{corollary:msplitn2zprime} to each of the two Appell functions found in the right-hand side of (\ref{equation:thetaLessOmega3}).  Furthermore, we will apply Corollary \ref{corollary:msplitn2zprime} to each Appell function twice:  once with $z^{\prime}=q^{8}$ and once with $z^{\prime}=q^{16}$.

We first specialize Corollary \ref{corollary:msplitn2zprime} with $(x;q)\to (q;q^6)$.  This yields
{\allowdisplaybreaks \begin{align}
m(q,z;q^{6})&=m(-q^{8},z';q^{24} )-q^{-5}m(-q^{-4},z';q^{24})
\label{equation:msplit2pre-z}\\
&\qquad +\frac{z'J_{12}^3}{j(qz;q^{6})j(z';q^{24})}
\Big [\frac{j(-q^{8}zz';q^{12})j(z^2/z';q^{24})}{j(-q^{8}z';q^{12})j(z;q^{12})}\notag\\
&\qquad \qquad -qz \frac{j(-q^{14}zz';q^{12})j(q^{12}z^2/z';q^{24})}{j(-q^{8}z';q^{12})j(q^{6}z;q^{12})}\Big ].\notag
\end{align}}%
To obtain the first Appell function on the right-hand side of (\ref{equation:thetaLessOmega3}), we further specialize (\ref{equation:msplit2pre-z}) to $z\to q^2$.  This gives
{\allowdisplaybreaks \begin{align}
m(q,q^{2};q^{6})&=m(-q^{8},z';q^{24} )-q^{-5}m(-q^{-4},z';q^{24})
\label{equation:omega3Appell-1}\\
&\qquad +\frac{z'J_{12}^3}{j(q^{3};q^{6})j(z';q^{24})}
\Big [\frac{j(-q^{10}z';q^{12})j(q^{4}/z';q^{24})}{j(-q^{8}z';q^{12})j(q^{2};q^{12})}\notag \\
&\qquad \qquad -q^{3} \frac{j(-q^{16}z';q^{12})j(q^{16}/z';q^{24})}{j(-q^{8}z';q^{12})j(q^{8};q^{12})}\Big ].\notag
\end{align}}%
To obtain the second Appell function on the right-hand side of (\ref{equation:thetaLessOmega3}), we further specialize (\ref{equation:msplit2pre-z}) to $z\to q^4$.  This gives
\begin{align}
m(q,q^{4};q^{6})&=m(-q^{8},z';q^{24} )-q^{-5}m(-q^{-4},z';q^{24})
\label{equation:omega3Appell-2}\\
&\qquad +\frac{z'J_{12}^3}{j(q^{5};q^{6})j(z';q^{24})}
\Big [\frac{j(-q^{12}z';q^{12})j(q^{8}/z';q^{24})}{j(-q^{8}z';q^{12})j(q^{4};q^{12})}\notag\\
&\qquad \qquad -q^{5} \frac{j(-q^{18}z';q^{12})j(q^{20}/z';q^{24})}{j(-q^{8}z';q^{12})j(q^{10};q^{12})}\Big ].\notag
\end{align}

%%%%
%%%% 
%%%%

For each of the two expansions (\ref{equation:omega3Appell-1}) and (\ref{equation:omega3Appell-2}), we will set $z^{\prime}=q^{8}$ and $z^{\prime}=q^{16}$.  This will give us four identities.  For the first expansion (\ref{equation:omega3Appell-1}), we first set $z^{\prime}\to q^{8}$ to get
\begin{align*}
m(q,q^{2};q^{6})&=m(-q^{8},q^{8};q^{24} )-q^{-5}m(-q^{-4},q^{8};q^{24})\\
&\qquad +\frac{q^{8}J_{12}^3}{j(q^{3};q^{6})j(q^{8};q^{24})}
\Big [\frac{j(-q^{18};q^{12})j(q^{-4};q^{24})}{j(-q^{16};q^{12})j(q^{2};q^{12})}
 -q^{3} \frac{j(-q^{24};q^{12})j(q^{8};q^{24})}{j(-q^{16};q^{12})j(q^{8};q^{12})}\Big ].
\end{align*}
We rewrite a theta function with (\ref{equation:1.7}) and then rewrite an Appell function with (\ref{equation:mxqz-flip}) and (\ref{equation:mxqz-fnq-z}) to get
\begin{align*}
m(q,q^{2};q^{6})&=m(-q^{8},q^{8};q^{24} )+q^{-1}m(-q^{4},q^{16};q^{24})\\
&\qquad +\frac{q^{8}J_{12}^3}{j(q^{3};q^{6})j(q^{8};q^{24})}
\Big [\frac{j(-q^{18};q^{12})j(q^{28};q^{24})}{j(-q^{16};q^{12})j(q^{2};q^{12})}
 -q^{3} \frac{j(-q^{24};q^{12})j(q^{8};q^{24})}{j(-q^{16};q^{12})j(q^{4};q^{12})}\Big ].
\end{align*}
Next, we rewrite several theta functions using (\ref{equation:1.8}).  This yields
\begin{align*}
m(q,q^{2};q^{6})&=m(-q^{8},q^{8};q^{24} )+q^{-1}m(-q^{4},q^{16};q^{24})\\
&\qquad -q^2\frac{J_{12}^3}{j(q^{3};q^{6})j(q^{8};q^{24})}
\Big [\frac{j(-q^{6};q^{12})j(q^{4};q^{24})}{j(-q^{4};q^{12})j(q^{2};q^{12})}
 +q \frac{j(-1;q^{12})j(q^{8};q^{24})}{j(-q^{4};q^{12})j(q^{4};q^{12})}\Big ].
\end{align*}
Then, we use (\ref{equation:1.12}) with $n=2$.  This gives
\begin{equation*}
j(q^{4};q^{24})=\frac{J_{24}}{J_{12}^2}j(q^{2};q^{12})j(-q^{2};q^{12}), \ 
j(q^{8};q^{24})=\frac{J_{24}}{J_{12}^2}j(q^{4};q^{12})j(-q^{4};q^{12}).
\end{equation*}
Cancelling out terms and pulling out common factors produces
\begin{align*}
m(q,q^{2};q^{6})&=m(-q^{8},q^{8};q^{24} )+q^{-1}m(-q^{4},q^{16};q^{24})\\
&\qquad -q^{2}\frac{J_{12}^3}{j(q^{3};q^{6})j(q^{8};q^{24})}\frac{J_{24}}{J_{12}^2j(-q^{4};q^{12})}\\
&\qquad \qquad \times \Big [j(-q^{6};q^{12})j(-q^{2};q^{12})
 +q j(-1;q^{12})j(-q^{4};q^{12})\Big ].
\end{align*}
Using (\ref{equation:1.7}) and (\ref{equation:H1Thm1.1}) with $(x,y;q)\to (-q,-q^5;q^6)$ yields
\begin{align*}
j(-q^{6};q^{12})&j(-q^{2};q^{12})
 +q j(-1;q^{12})j(-q^{4};q^{12})\\
 &=j(-q^{6};q^{12})j(-q^{10};q^{12})  +q j(-q^{12};q^{12})j(-q^{4};q^{12})
 =j(-q;q^{6})j(-q^5;q^6).
\end{align*}
Substituting in gives
\begin{align*}
m(q,q^{2};q^{6})&=m(-q^{8},q^{8};q^{24} )+q^{-1}m(-q^{4},q^{16};q^{24})\\
&\qquad -q^{2}\frac{J_{12}^3}{j(q^{3};q^{6})j(q^{8};q^{24})}\frac{J_{24}}{J_{12}^2j(-q^{4};q^{12})}
\cdot j(-q;q^{6})j(-q^5;q^6).
\end{align*}
In more compact notation, the first of our four identities then reads
\begin{equation}
m(q,q^{2};q^{6})=m(-q^{8},q^{8};q^{24} )+q^{-1}m(-q^{4},q^{16};q^{24})
 -q^{2}\frac{J_{12}^3}{J_{3,6}J_{8}}\frac{J_{24}}{J_{12}^2\overline{J}_{4,12}}
\cdot \overline{J}_{1,6}^2.
\label{equation:Omega-expansion1}
\end{equation}

%%%%
%%%%
%%%%

For the first expansion (\ref{equation:omega3Appell-1}), we now set $z^{\prime}=q^{16}$ to get
\begin{align*}
m(q,q^{2};q^{6})&=m(-q^{8},q^{16};q^{24} )-q^{-5}m(-q^{-4},q^{16};q^{24})\\
&\qquad +\frac{q^{16}J_{12}^3}{j(q^{3};q^{6})j(q^{16};q^{24})}
\Big [\frac{j(-q^{26};q^{12})j(q^{-12};q^{24})}{j(-q^{24};q^{12})j(q^{2};q^{12})}
 -q^{3} \frac{j(-q^{32};q^{12})j(1;q^{24})}{j(-q^{24};q^{12})j(q^{8};q^{12})}\Big ].
\end{align*}
The second quotient of theta functions vanishes.  We rewrite the second Appell function using (\ref{equation:mxqz-flip}) and (\ref{equation:mxqz-fnq-z}).  We then use (\ref{equation:1.7}) to rewrite a theta function in the numerator of the remaining theta quotient to get
\begin{align*}
m(q,q^{2};q^{6})&=m(-q^{8},q^{16};q^{24} )+q^{-1}m(-q^{4},q^{8};q^{24})\\
&\qquad +\frac{q^{16}J_{12}^3}{j(q^{3};q^{6})j(q^{8};q^{24})}
\cdot \frac{j(-q^{26};q^{12})j(q^{36};q^{24})}{j(-q^{24};q^{12})j(q^{2};q^{12})}.
\end{align*}
We rewrite several theta functions by using (\ref{equation:1.8}) to get
\begin{align*}
m(q,q^{2};q^{6})&=m(-q^{8},q^{16};q^{24} )+q^{-1}m(-q^{4},q^{8};q^{24})\\
&\qquad -\frac{J_{12}^3}{j(q^{3};q^{6})j(q^{8};q^{24})}
\cdot  \frac{j(-q^{2};q^{12})j(q^{12};q^{24})}{j(-1;q^{12})j(q^{2};q^{12})}.
\end{align*}
In more compact notation, the second of our four identities then reads
\begin{equation}
m(q,q^{2};q^{6})=m(-q^{8},q^{16};q^{24} )+q^{-1}m(-q^{4},q^{8};q^{24})
-\frac{J_{12}^3}{J_{3,6}J_{8}}
\cdot  \frac{\overline{J}_{2,12}J_{12,24}}{\overline{J}_{0,12}J_{2,12}}.
\label{equation:Omega-expansion2}
\end{equation}

%%%%
%%%%
For the second expansion (\ref{equation:omega3Appell-2}), we first set $z^{\prime}=q^{8}$ to get
\begin{align*}
m(q,q^{4};q^{6})&=m(-q^{8},q^{8};q^{24} )-q^{-5}m(-q^{-4},q^{8};q^{24})\\
&\qquad +\frac{q^{8}J_{12}^3}{j(q^{5};q^{6})j(q^{8};q^{24})}
\Big [\frac{j(-q^{20};q^{12})j(1;q^{24})}{j(-q^{16};q^{12})j(q^{4};q^{12})}
 -q^{5} \frac{j(-q^{26};q^{12})j(q^{12};q^{24})}{j(-q^{16};q^{12})j(q^{10};q^{12})}\Big ].
\end{align*}
The first quotient of theta functions vanishes.  We rewrite the second Appell function using (\ref{equation:mxqz-flip}) and (\ref{equation:mxqz-fnq-z}).  Next, we rewrite a theta function in the denominator using (\ref{equation:1.7}) to get
\begin{align*}
m(q,q^{4};q^{6})&=m(-q^{8},q^{8};q^{24} )+q^{-1}m(-q^{4},q^{16};q^{24})\\
&\qquad -\frac{q^{13}J_{12}^3}{j(q;q^{6})j(q^{8};q^{24})}
\cdot \frac{j(-q^{26};q^{12})j(q^{12};q^{24})}{j(-q^{16};q^{12})j(q^{2};q^{12})}.
\end{align*}
We then apply (\ref{equation:1.8}) to several theta functions to get
\begin{align*}
m(q,q^{4};q^{6})&=m(-q^{8},q^{8};q^{24} )+q^{-1}m(-q^{4},q^{16};q^{24})\\
&\qquad -q\frac{J_{12}^3}{j(q;q^{6})j(q^{8};q^{24})}
\cdot \frac{j(-q^{2};q^{12})j(q^{12};q^{24})}{j(-q^{4};q^{12})j(q^{2};q^{12})}.
\end{align*}
In more compact notation, the third of our four identities then reads
\begin{equation}
m(q,q^{4};q^{6})=m(-q^{8},q^{8};q^{24} )+q^{-1}m(-q^{4},q^{16};q^{24})\\
 -q\frac{J_{12}^3}{J_{1,6}J_{8}}
\cdot \frac{\overline{J}_{2,12}J_{12,24}}{\overline{J}_{4,12}J_{2,12}}.
\label{equation:Omega-expansion3}
\end{equation}

%%%%
%%%%
%%%%
For the second expansion (\ref{equation:omega3Appell-2}), we next set $z^{\prime}=q^{16}$ to get
\begin{align*}
m(q,q^{4};q^{6})&=m(-q^{8},q^{16};q^{24} )-q^{-5}m(-q^{-4},q^{16};q^{24})\\
&\qquad +\frac{q^{16}J_{12}^3}{j(q^{5};q^{6})j(q^{16};q^{24})}
\Big [\frac{j(-q^{28};q^{12})j(q^{-8};q^{24})}{j(-q^{24};q^{12})j(q^{4};q^{12})}
 -q^{5} \frac{j(-q^{34};q^{12})j(q^{4};q^{24})}{j(-q^{24};q^{12})j(q^{10};q^{12})}\Big ].
\end{align*}
We rewrite the second Appell function using (\ref{equation:mxqz-flip}) and (\ref{equation:mxqz-fnq-z}).  We then use (\ref{equation:1.7}) to rewrite a few theta functions.  This yields
\begin{align*}
m(q,q^{4};q^{6})&=m(-q^{8},q^{16};q^{24} )+q^{-1}m(-q^{4},q^{8};q^{24})\\
&\qquad +\frac{q^{16}J_{12}^3}{j(q^{5};q^{6})j(q^{8};q^{24})}
\Big [\frac{j(-q^{28};q^{12})j(q^{32};q^{24})}{j(-q^{24};q^{12})j(q^{4};q^{12})}
 -q^{5} \frac{j(-q^{34};q^{12})j(q^{4};q^{24})}{j(-q^{24};q^{12})j(q^{2};q^{12})}\Big ].
\end{align*}
Applying (\ref{equation:1.8}) to several theta functions produces
\begin{align*}
m(q,q^{4};q^{6})&=m(-q^{8},q^{16};q^{24} )+q^{-1}m(-q^{4},q^{8};q^{24})\\
&\qquad -\frac{J_{12}^3}{j(q^{5};q^{6})j(q^{8};q^{24})}
\Big [\frac{j(-q^{4};q^{12})j(q^{8};q^{24})}{j(-1;q^{12})j(q^{4};q^{12})}
 +q\frac{j(-q^{10};q^{12})j(q^{4};q^{24})}{j(-1;q^{12})j(q^{2};q^{12})}\Big ].
\end{align*}
Again, we use (\ref{equation:1.12}) with $n=2$.  This gives
\begin{equation*}
j(q^{8};q^{24})=\frac{J_{24}}{J_{12}^2}j(q^{4};q^{12})j(-q^{4};q^{12}), \
j(q^{4};q^{24})=\frac{J_{24}}{J_{12}^2}j(q^{2};q^{12})j(-q^{2};q^{12}).
\end{equation*}
Substituting in and pulling out common factors brings us to
\begin{align*}
m(q,q^{4};q^{6})&=m(-q^{8},q^{16};q^{24} )+q^{-1}m(-q^{4},q^{8};q^{24})\\
&\qquad -\frac{J_{12}^3}{j(q^{5};q^{6})j(q^{8};q^{24})} \frac{J_{24}}{J_{12}^2j(-1;q^{12})}\\
&\qquad \qquad \times\Big [j(-q^{4};q^{12})j(-q^{4};q^{12})
 +qj(-q^{10};q^{12})j(-q^{2};q^{12})\Big ].
\end{align*}
Using (\ref{equation:1.7}) and (\ref{equation:H1Thm1.1})  with $(x,y;q)\to (-q,-q^3;q^6)$ yields
\begin{align*}
j(-q^{4};q^{12})&j(-q^{4};q^{12}) +qj(-q^{10};q^{12})j(-q^{2};q^{12})\\
&=j(-q^{4};q^{12})j(-q^{4};q^{12}) +qj(-q^{10};q^{12})j(-q^{10};q^{12})
=j(-q;q^{6})j(-q^3;q^6).
\end{align*}
Substituting in gives
\begin{align*}
m(q,q^{4};q^{6})&=m(-q^{8},q^{16};q^{24} )+q^{-1}m(-q^{4},q^{8};q^{24})\\
&\qquad -\frac{J_{12}^3}{j(q^{5};q^{6})j(q^{8};q^{24})} \frac{J_{24}}{J_{12}^2j(-1;q^{12})}
 \cdot j(-q;q^{6})j(-q^3;q^6).
\end{align*}
In more compact notation, the fourth of our four identities then reads
\begin{equation}
m(q,q^{4};q^{6})=m(-q^{8},q^{16};q^{24} )+q^{-1}m(-q^{4},q^{8};q^{24})
 -\frac{J_{12}^3}{J_{1,6}J_{8}} \frac{J_{24}}{J_{12}^2\overline{J}_{0,12}}
 \cdot \overline{J}_{1,6}\overline{J}_{3,6}.
 \label{equation:Omega-expansion4}
\end{equation}

%%%%%
%%%%%
%%%%%

We add the four identities (\ref{equation:Omega-expansion1}), (\ref{equation:Omega-expansion2}), (\ref{equation:Omega-expansion3}), (\ref{equation:Omega-expansion4}) to get
{\allowdisplaybreaks \begin{align*}
2m(q,q^{2};q^{6})+2m(q,q^{4};q^{6})
&=2m(-q^{8},q^{8};q^{24} )+2q^{-1}m(-q^{4},q^{16};q^{24})
\\
%%%
&\quad +2m(-q^{8},q^{16};q^{24} )+2q^{-1}m(-q^{4},q^{8};q^{24})
\\
%%%
&\quad 
 -\frac{J_{12}^3}{J_{3,6}J_{8}}
\cdot  \frac{\overline{J}_{2,12}J_{12,24}}{\overline{J}_{0,12}J_{2,12}}-q\frac{J_{12}^3}{J_{1,6}J_{8}}
\cdot \frac{\overline{J}_{2,12}J_{12,24}}{\overline{J}_{4,12}J_{2,12}}\\
%%%
&\quad 
 -\frac{J_{12}^3}{J_{1,6}J_{8}} \frac{J_{24}}{J_{12}^2\overline{J}_{0,12}}
 \cdot \overline{J}_{1,6}\overline{J}_{3,6}
  -q^{2}\frac{J_{12}^3}{J_{3,6}J_{8}}\frac{J_{24}}{J_{12}^2\overline{J}_{4,12}}
\cdot \overline{J}_{1,6}^2.
\end{align*}}%
Using the Appell function expressions for the two third order mock theta functions $f_{3}(q)$ and $\omega_{3}(q)$ as found in (\ref{equation:thetaLessf3}) and (\ref{equation:thetaLessOmega3}), we can rewrite the above as
{\allowdisplaybreaks \begin{align*}
-2q\omega_{3}(q)
&=f_{3}(q^8)+2q^{3}\omega_{3}(-q^4) 
 -\frac{J_{12}^3}{J_{3,6}J_{8}}
\cdot  \frac{\overline{J}_{2,12}J_{12,24}}{\overline{J}_{0,12}J_{2,12}}-q\frac{J_{12}^3}{J_{1,6}J_{8}}
\cdot \frac{\overline{J}_{2,12}J_{12,24}}{\overline{J}_{4,12}J_{2,12}}\\
%%%
&\quad 
 -\frac{J_{12}^3}{J_{1,6}J_{8}} \frac{J_{24}}{J_{12}^2\overline{J}_{0,12}}
 \cdot \overline{J}_{1,6}\overline{J}_{3,6}
  -q^{2}\frac{J_{12}^3}{J_{3,6}J_{8}}\frac{J_{24}}{J_{12}^2\overline{J}_{4,12}}
\cdot \overline{J}_{1,6}^2.
\end{align*}}%

Let us combine the theta quotients pairwise.  We combine fractions and then use (\ref{equation:1.12}) with $n=2$ to combine two pairs of theta functions.  This yields
{\allowdisplaybreaks \begin{align*}
\frac{J_{12}^3}{J_{1,6}J_{8}} \frac{J_{24}}{J_{12}^2\overline{J}_{0,12}}
 \cdot \overline{J}_{1,6}\overline{J}_{3,6}
& +q^{2}\frac{J_{12}^3}{J_{3,6}J_{8}}\frac{J_{24}}{J_{12}^2\overline{J}_{4,12}}
\cdot \overline{J}_{1,6}^2\\
&=\frac{J_{12}^3}{J_{8}} \frac{J_{24}}{J_{12}^2}
 \cdot \overline{J}_{1,6}
 \left ( \frac{\overline{J}_{3,6}}{J_{1,6} \overline{J}_{0,12}}
 +q^2 \frac{\overline{J}_{1,6}}{J_{3,6} \overline{J}_{4,12}} \right )\\
&=\frac{J_{12}^3}{J_{8}} \frac{J_{24}}{J_{12}^2}
 \cdot \overline{J}_{1,6}
 \left ( \frac{\overline{J}_{3,6}J_{3,6}\overline{J}_{4,12} +q^2\overline{J}_{1,6}J_{1,6} \overline{J}_{0,12}}
 {J_{1,6} J_{3,6}\overline{J}_{0,12}\overline{J}_{4,12}}
\right )\\
&=\frac{J_{12}^3}{J_{8}} \frac{J_{24}}{J_{12}^2}
 \cdot \overline{J}_{1,6}\frac{J_{6}^2}{J_{12}}
 \left ( \frac{J_{6,12}\overline{J}_{4,12} +q^2J_{2,12} \overline{J}_{0,12}}
 {J_{1,6} J_{3,6}\overline{J}_{0,12}\overline{J}_{4,12}}
\right ).
\end{align*}}%
Using (\ref{equation:H1Thm1.1}) with $(x,y;q)\to (-q^2,-q^2;-q^6)$ produces
\begin{align*}
j(-q^2;-q^6)j(-q^2;-q^6)=j(-q^4;q^{12})j(q^{6};q^{12})+q^2j(q^{10};q^{12})j(-1,q^{12}).
\end{align*}
Hence
\begin{align*}
\frac{J_{12}^3}{J_{1,6}J_{8}} \frac{J_{24}}{J_{12}^2\overline{J}_{0,12}}
 \cdot \overline{J}_{1,6}\overline{J}_{3,6}
& +q^{2}\frac{J_{12}^3}{J_{3,6}J_{8}}\frac{J_{24}}{J_{12}^2\overline{J}_{4,12}}
\cdot \overline{J}_{1,6}^2\\
&=\frac{J_{12}^3}{J_{8}} \frac{J_{24}}{J_{12}^2}
 \cdot \overline{J}_{1,6}\frac{J_{6}^2}{J_{12}}
 \left ( \frac{J_{4}j(-q^2;q^4)}
 {J_{1,6} J_{3,6}\overline{J}_{0,12}\overline{J}_{4,12}}
\right ).
\end{align*}
Product rearrangements give
\begin{equation}
\frac{J_{12}^3}{J_{1,6}J_{8}} \frac{J_{24}}{J_{12}^2\overline{J}_{0,12}}
 \cdot \overline{J}_{1,6}\overline{J}_{3,6}
  +q^{2}\frac{J_{12}^3}{J_{3,6}J_{8}}\frac{J_{24}}{J_{12}^2\overline{J}_{4,12}}
\cdot \overline{J}_{1,6}^2
=\frac{1}{2}\frac{J_{2}J_{4}^6}{J_{1}^2J_{8}^4}.
\label{equation:Watson-firstThetaPair}
\end{equation}
%%%%%
%%%%%
%%%%%

Let us consider the second pair of theta quotients.  We combine fractions and then use (\ref{equation:1.10n2}) to break up two theta functions.  This gives
{\allowdisplaybreaks \begin{align*}
\frac{J_{12}^3}{J_{3,6}J_{8}}
\cdot  \frac{\overline{J}_{2,12}J_{12,24}}{\overline{J}_{0,12}J_{2,12}}
&+q\frac{J_{12}^3}{J_{1,6}J_{8}}
\cdot \frac{\overline{J}_{2,12}J_{12,24}}{\overline{J}_{4,12}J_{2,12}}\\
&=\frac{J_{12}^3}{J_{8}}
\cdot  \frac{\overline{J}_{2,12}J_{12,24}}{J_{2,12}}
\left (\frac{1}{J_{3,6}\overline{J}_{0,12}}+q\frac{1}{J_{1,6} \overline{J}_{4,12}} \right ) \\
&=\frac{J_{12}^3}{J_{8}}
\cdot  \frac{\overline{J}_{2,12}J_{12,24}}{J_{2,12}}
\left (\frac{J_{1,6} \overline{J}_{4,12}+qJ_{3,6}\overline{J}_{0,12}}
{J_{3,6}\overline{J}_{0,12}J_{1,6} \overline{J}_{4,12}} \right )\\
&=\frac{J_{12}^3}{J_{8}}
\cdot  \frac{\overline{J}_{2,12}J_{12,24}}{J_{2,12}}
\frac{J_{6}}{J_{12}^2}\left (\frac{J_{1,12} J_{7,12}\overline{J}_{4,12}+qJ_{3,12}J_{9,12}\overline{J}_{0,12}}
{J_{3,6}\overline{J}_{0,12}J_{1,6} \overline{J}_{4,12}}\right ).
\end{align*}}%
We use the Weierstrass relation (\ref{equation:3termWeier}) with $(a,b,c,d;q)\to (q^6,q^3,q^2,-q^2;q^{12})$ to get
\begin{align*}
&j(q^{8};q^{12})j(q^{4};q^{12})j(-q^{5};q^{12})j(-q;q^{12})\\
&=j(-q^{8};q^{12})j(-q^{4};q^{12})j(q^{5};q^{12})j(q;q^{12})
+qj(q^{9};q^{12})j(q^{3};q^{12})j(-q^{4};q^{12})j(-1;q^{12}).
\end{align*}
Hence
{\allowdisplaybreaks \begin{align*}
\frac{J_{12}^3}{J_{3,6}J_{8}}
\cdot  \frac{\overline{J}_{2,12}J_{12,24}}{\overline{J}_{0,12}J_{2,12}}
&+q\frac{J_{12}^3}{J_{1,6}J_{8}}
\cdot \frac{\overline{J}_{2,12}J_{12,24}}{\overline{J}_{4,12}J_{2,12}}\\
&=\frac{J_{12}^3}{J_{8}}
\cdot  \frac{\overline{J}_{2,12}J_{12,24}}{J_{2,12}}
\frac{J_{6}}{J_{12}^2}\left (\frac{J_{4}^2\overline{J}_{5,12}\overline{J}_{1,12}}
{J_{3,6}\overline{J}_{0,12}J_{1,6} \overline{J}_{4,12}\overline{J}_{4,12}} \right ).
\end{align*}}%
Product rearrangements give
\begin{equation*}
\frac{J_{12}^3}{J_{3,6}J_{8}}
\cdot  \frac{\overline{J}_{2,12}J_{12,24}}{\overline{J}_{0,12}J_{2,12}}
+q\frac{J_{12}^3}{J_{1,6}J_{8}}
\cdot \frac{\overline{J}_{2,12}J_{12,24}}{\overline{J}_{4,12}J_{2,12}}=\frac{1}{2}\frac{J_{2}J_{4}^6}{J_{1}^2J_{8}^4}.
\label{equation:Watson-secondThetaPair}
\end{equation*}
Hence we can write
{\allowdisplaybreaks \begin{align*}
-2q\omega_{3}(q)
&=f_{3}(q^8)+2q^{3}\omega_{3}(-q^4)-\frac{J_{2}J_{4}^6}{J_{1}^2J_{8}^4},
\end{align*}}%
and the result follows.
%%%%%%%%
%%%%%%%%
%%%%%%%%

\section*{Acknowledgements}
The work is supported by the Ministry of Science and Higher Education of the Russian
Federation (agreement no. 075-15-2025-343).

\end{document}